\documentclass[a4,10pt]{article}
\usepackage[a4paper,left=2cm,right=2cm,top=2cm,bottom=2cm]{geometry}
\usepackage{amsmath,amssymb}
\usepackage{hyperref}
\usepackage{setspace}
\usepackage[utf8]{inputenc}
\usepackage[english]{babel}
\usepackage{amsfonts, amstext, amsthm, booktabs, dcolumn}
\usepackage{lscape}
\usepackage[T1]{fontenc}
\usepackage{authblk,color}
\usepackage{booktabs,multirow}

\newtheorem{proposition}{\bf Proposition}[section]
\newtheorem{theorem}{\bf Theorem}[section]

\newtheorem{corollary}{\bf Corollary}[section]

\newtheorem{lemma}{\bf Lemma}[section]

\newcommand\floor[1]{\left\lfloor#1\right\rfloor}

\begin{document}

\title{Pseudo-binomial Approximation to $(k_1, k_2)$-runs}
\author[*]{A. N. Kumar}
\author[**]{N. S. Upadhye}
\affil[ ]{\small Department of Mathematics, Indian Institute of Technology Madras,}
\affil[ ]{\small Chennai-600036, India.}
\affil[*]{\small Email: amit.kumar2703@gmail.com}
\affil[**]{\small Email: neelesh@iitm.ac.in}
\date{}
\maketitle

\begin{abstract}
\noindent
($k_1,k_2$)-runs have received a special attention in the literature and its distribution can be obtained using combinatorial method (Huang and Tsai \cite{HT}) and Markov chain approach (Dafnis {\em et al.} \cite{DAP}). But the formulae are difficult to use when the number of Bernoulli trials is too large under identical setup and is generally intractable under non-identical setup. So, it is useful to approximate it with a suitable random variable. In this paper, it is demonstrated that pseudo-binomial is most suitable distribution for approximation and the approximation results are derived using Stein's method. Also, application of these results is demonstrated through real-life problems. It is shown that the bounds obtained are either comparable to or improvement over bounds available in the literature.
\end{abstract}

\noindent
\maketitle

\section{Introduction and Preliminaries}\label{two:intro}
Let $\xi_1, \xi_2, \dotsc,\xi_n$ be a sequence of independent Bernoulli trials with success probability ${\mathbb P}(\xi_l=1)=p_l$, $1 \le l \le n$. A $(k_1,k_2)$-event is a pattern of at least $k_1$ consecutive failures followed by at least $k_2$ consecutive successes, where $(k_1,k_2)$ is a pair of non-negative integers, including $0$, excluding $(0,0)$. Let $B_{k_1,k_2}^n$ be the number of occurrences of $(k_1,k_2)$-events in $n$ trials, then it is called a modified binomial distribution of order $k$ (see, Huang and Tsai \cite{HT}) or $(k_1,k_2)$-runs (see, Upadhye {\em et al.} \cite{UCV} and Vellaisamy \cite{V}). Mathematically, $(k_1,k_2)$-runs can be represented as a dependent setup that arises from an independent sequence of Bernoulli trials as follows:
\begin{equation}
I_l = (1-\xi_{l})\dotsb (1-\xi_{l+k_1-1}) \xi_{l+k_1}\dotsb \xi_{l+k_1+k_2-1} \quad{\rm and}\quad B_{k_1,k_2}^n = \sum_{l=1}^{n-k_1-k_2+1} I_l.\label{two:ij}
\end{equation}
The distributional properties of $B_{k_1,k_2}^n$, such as probability generating function (PGF), probability mass function (PMF), waiting time distribution and their moments are studied by Huang and Tsai \cite{HT} and Dafnis {\em et al.} \cite{DAP} for independent and identically distributed (iid) trials. Also, the exact distribution of ($k_1,k_2$)-runs can be obtained using Markov chain approach given by Fu and Koutras \cite{FK} for identical and non-identical trials. But the formulae (or recursive relations) are not practically useful whenever $n$, $k_1$, and $k_2$ are too large for identical trials and the distribution of $B_{k_1,k_2}^n$ is generally intractable for non-identical trials. It is of interest to study, the approximation problem related to $B_{k_1,k_2}^n$.\\
Approximations to runs are widely studied in the literature, for example, negative binomial approximation to $k$-runs (Wang and Xia \cite{WX}), Poisson approximation to $(k_1,k_2)$-events (Vellaisamy \cite{V}) and Poisson approximations for the reliability of the system (Godbole \cite{G}). Recently, Kumar and Upadhye \cite{KU} obtained bounds between negative binomial and a function of waiting time for $(k_1, k_2)$-events and Upadhye {\em et al.} \cite{UCV} derived a bound for binomial convoluted Poisson approximation to $(1,1)$-runs. For more details and applications of runs, see Aki {\em et al.} \cite{AKH}, Antzoulakos {\em et al.} \cite{ABK}, Ankzoulakos and Chadjiconstantinidis \cite{AC}, Balakrishnan and Koutras \cite{B}, Makri {\em et al.} \cite{MPP}, Philippou {\em et al.} \cite{PGP}, Philippou and Makri \cite{PM} and references therein.\\
Next, we discuss Stein's method (Stein \cite{stein}) which is an important tool to obtain the results derived in this paper. As a first step, we identify a suitable operator (called Stein operator, say ${\cal A}_X$) for a random variable $X$ which is acting on a class of functions (say ${\cal G}_X$) such that ${\mathbb E} \left[{\cal A}_X g(X)\right]=0,~{\rm for}~g \in {\cal G}_X.$\\
The next step is to obtain the solution of the following equation
$${\cal A}_X g(m) = f (m) - {\mathbb E}f (X),\quad m \in {\mathbb Z} ~{\rm and}~ f \in {\cal G}_X,$$
which is known as Stein equation. Finally, substituting a random variable $\tilde{X}$ for $m$ in Stein equation, and taking expectations and supremum, we get
$$d_{TV}\left(X,\tilde{X}\right) := \sup_{f \in {\cal H}}\left|{\mathbb E}f(\tilde{X}) - {\mathbb E}f(X)\right| = \sup_{f \in {\cal H}}\left|{\mathbb E}\left[{\cal A }_Xg(\tilde{X})\right]\right|,$$
where ${\cal H} = \{{\bf 1}(S)|~ S~measurable\}$ and ${\bf 1}(S)$ is the indicator function of the set $S$. Throughout this paper, let $\cal G$ be the set of all bounded functions and
\begin{equation}
{\cal G}_X = \left\{g|~g \in {\cal G}~{\rm such~that}~g(0)=0~{\rm and} ~g(x)=0, ~{\rm for}~x \notin {Supp(X)}  \right\}\label{two:gx}
\end{equation}
be associated with Stein operator ${\cal A}_X$, where $Supp(X)$ denotes the support of a random variable $X$.  For more details and applications, see Barbour {\em et al.} \cite{BCL,BHJ}, Barbour and Chen \cite{BAC}, Chen {\em et al.} \cite{CGS}, \v{C}ekanavi\v{c}ius \cite{CV}, Norudin and Peccati \cite{NP}, Reinert \cite{R1} and references therein. A Stein operator can be obtained using different approaches available in the literature (see Stein \cite{stein}, Barbour and G\"{o}tze \cite{Ba,Go}, Diaconis and Zabell \cite{DZ}, Ley {\em et al.} \cite{LRS} and Upadhye {\em et al.} \cite{UCV}). However, for identical trials, we focus on PGF approach given by Upadhye {\em et al.} \cite{UCV} and derive a Stein operator for $B_{k_1,k_2}^n$ as a perturbation of pseudo-binomial operator using the technique discussed by Kumar and Upadhye \cite{KU}. For non-identical trials, PGF for such distribution is generally complicated. So, PGF approach becomes difficult to apply and involves unnecessary complications. Hence, we use coupling approach (see Wang and Xia \cite{WX}) to obtain error bounds for non-identical trials. It is worth to note that PGF approach uses more information about the distribution and gives better approximation results. \\
Next, from \eqref{two:ij}, it can be verified that
\begin{equation}
{\mathbb E}\big(B_{k_1,k_2}^n\big)=\sum_{l=1}^{n-k_1-k_2+1}\hspace{-0.3cm}a(p_l)\quad {\rm and}\quad Var\big(B_{k_1,k_2}^n\big)=\sum_{l=1}^{n-k_1-k_2+1}\hspace{-0.3cm}a(p_l)-\sum_{l=1}^{n-k_1-k_2+1}\hspace{-0.3cm}(a(p_l))^2-2\sum_{\substack{l<r\\r-l\le k_1+k_2-1}}\hspace{-0.5cm}a(p_l)a(p_r),\label{two:MVAR}
\end{equation}
where $a(p_l):={\mathbb E}(I_l)=(1-p_l)\dotsb(1-p_{l+k_1-1})p_{l+k_1}\dotsb p_{l+k_1+k_2-1}$. From \eqref{two:MVAR}, it is clear that ${\mathbb E}\big(B_{k_1,k_2}^n\big)>Var\big(B_{k_1,k_2}^n\big)$ and hence pseudo-binomial is suitable to do approximation for the distribution of $B_{k_1,k_2}^n$. Now, let $Z$ follow pseudo-binomial distribution (see Upadhye {\em et al.} \cite{UCV}) with PMF
\begin{equation}
{\mathbb P}(Z=m)= \frac{1}{\cal C} {\alpha \choose m} \check{p}^m \check{q}^{\alpha-m}, \quad m=0,1,2,\dotsc,\floor{\alpha },\label{two:ppmf}
\end{equation}
where $\alpha > 0$ and $0 < \check{p}<1$ with $\check{q}=1-\check{p}$ be the parameter of $Z$, $\floor{\alpha}$ is the greatest integer function of $\alpha$ and ${\cal C} = \sum_{m=0}^{\floor{\alpha}}{\alpha \choose m} \check{p}^m \check{q}^{\alpha-m}$. From ($5$) of Upadhye {\em et al.} \cite{UCV}, Stein operator for $Z$ is given by
\begin{equation}
{\cal A}_0g(m) = (\alpha-m)\check{p} g(m+1) - m \check{q} g(m), \quad m=0, 1, \dotsc, \floor{\alpha}\label{two:sboper}
\end{equation}
\vspace{-0.19cm}
and the bound for the solution of the Stein equation is
\begin{equation}
\|\Delta g\| \le \frac{2 \|f\|}{\floor{\alpha} \check{p}\check{q}}.\label{two:bound}
\end{equation}
The paper is organized as follows. In Section \ref{two:ARiid}, for identical trials, we first obtain recursive relations in PGF. Next, a Stein operator for $B_{k_1,k_2}^n$ is obtained via PGF approach. Furthermore,  a Stein operator for $B_{k_1,k_2}^n$ is shown as a perturbation of pseudo-binomial operator. Finally, we obtain the approximation results by comparing Stein operators. In Section \ref{two:ARit}, for non-identical trials, we first obtain the one-parameter approximation result. Next, using an appropriate coupling, we derive error bounds for two-parameter approximation.  In Section \ref{two:CR}, we compare our results to the existing results  and make some relevant remarks. In Section \ref{two:app}, we give applications of the results and compare bounds with the known bounds in the literature.

\section{Approximation Results for iid Trials}\label{two:ARiid}
In this section, we first obtain recursive relations in PGF and its derivative for the random variable $B_{k_1,k_2}^n$. During this process, an explicit form of PMF is also obtained by solving the recursive relation given by Dafnis {\em et al.} \cite{DAP}. Next, we derive a Stein operator for $B_{k_1,k_2}^n$ via PGF approach. Also, a Stein operator which can be seen as a perturbation of pseudo-binomial operator is obtained. Finally, using perturbation technique, error bounds between pseudo-binomial and $B_{k_1,k_2}^n$ are derived by matching first and first two moments.\\
Throughout in this section, let $p_{i}=p,$ for all $i=1,2,\dotsc,n$,
$$k:=k_1+k_2\quad{\rm and}\quad a(p_l)=(1-p)^{k_1}p^{k_1}=:a(p), ~{\rm for~all}~l=1,2,\dotsc,n-k+1.$$

\subsection{Recursive Relations}
Let $\phi_n(t)$ denote the PGF of $B_{k_1,k_2}^n$. Then, from Theorem $1$ of Huang and Tsai \cite{HT}, we have
$$\phi_n(t):=\phi_n(t;~k_1,~k_2)=\sum_{m=0}^{\floor{n/k}}p_{m,n} t^m=\sum_{i=i_0}^{n}{i \choose \alpha_i}\left[a(p)(t-1)\right]^{i-\alpha_i}$$
or equivalently
\begin{equation}
\phi_n(t)  = \sum_{m=0}^{\floor{n/k}}{{n-m(k-1)} \choose m} (a(p)(t-1))^m,\label{two:pgf}
\end{equation}
where $p_{m, n} := {\mathbb P}\big(B_{k_1,k_2}^n=m\big)$ is the PMF of $B_{k_1,k_2}^n$, $i_0=-\floor{-n/k}$ and $\alpha_i=(k i -n)/(k-1)$ if the right-hand side is an integer, otherwise $-1$ with ${b \choose -1}=0$ for any $b$. Also, let $\phi(z,t)$ denote the double PGF of $B_{k_1,k_2}^n$ (see Dafnis {\em et al.} \cite{DAP} and Huang and Tsai \cite{HT} for details) and is given by
\begin{equation}
\phi(z,t):= \phi(z;~t,~k_1,~k_2) = \sum_{n=0}^{\infty} \phi_n(t) z^n = \frac{1}{1-z-a(p)(t-1)z^k}.\label{two:dpgf}
\end{equation}

\begin{lemma}\label{two:relpgf}
Let $n \ge k$, then $\phi_n(t)$ satisfies the following relations:
\begin{itemize}
\item[{\rm (i)}] $\phi_n^\prime(t) = (n-k+1) a(p) \phi_{n-k}(t) - a(p)(k-1)(t-1)\phi_{n-k}^\prime(t)$.
\item[{\rm (ii)}] $\phi_{n-1}^{\prime}(t) = (n-k) a(p) \phi_{n-k}(t) - a(p) k (t-1) \phi_{n-k}^{\prime}(t)$.
\item[{\rm (iii)}] $\phi_n^\prime(t) = \phi_{n-1}^\prime(t) + a(p) \phi_{n-k}(t) + a(p)(t-1)\phi_{n-k}^\prime(t).$
\item[{\rm (iv)}] $\displaystyle{\left[1+\tilde{p}(t-1)\right]\phi_{n}^{\prime}(t) =\frac{n}{k} \tilde{p} \phi_n(t) - a(p) \sum_{u=0}^{k-2} \left(\frac{n+u+1}{k-1}\right) \left(\frac{k}{k-1}\right)^{u}\phi_{n-k+u+1}(t)}$,\\
where $\displaystyle{\tilde{p} = k a(p)\left(\frac{k}{k-1}\right)^{k-1} }$ and $\phi_n^\prime$ denotes the first derivative of $\phi_n$  w.r.t. $t$.
\end{itemize}
\end{lemma}

\noindent
{\bf Proof}. (i) Differentiating (\ref{two:pgf}) w.r.t. $t$, we have
\begin{align}
\phi_n^\prime(t) &= \sum_{m=1}^{\floor{n/k}} m {{n-m(k-1)} \choose m} a(p)^m (t-1)^{m-1}\label{two:pgf1}~~~~~~~~~~~~~~~~~~~~~~~~~~~~~~~~~~~~~~~~~~~~~~~
\end{align}
\begin{align}
&=  \sum_{m=1}^{\floor{n/k}} (n-m(k-1)) {{n-m(k-1)-1} \choose m-1} a(p)^m (t-1)^{m-1}\nonumber\\
&= \sum_{m=0}^{\floor{n/k}-1} (n-(m+1)(k-1)) {{n-(m+1)(k-1)-1} \choose m} a(p)^{m+1} (t-1)^{m}\nonumber\\
&=(n-k+1)a(p)\sum_{m=0}^{\floor{(n-k)/k}}{n-k-m(k-1)\choose m}(a(p)(t-1))^m\nonumber\\
&~~~-a(p)(k-1)(t-1)\sum_{m=1}^{\floor{(n-k)/k}}m{n-k-m(k-1)\choose m}(a(p))^mt^{m-1}.\nonumber
\end{align}
Using \eqref{two:pgf} and \eqref{two:pgf1} with $n=n-k$, the result follows.\\
(ii) Differentiating \eqref{two:dpgf} w.r.t $z$ and $t$, we have
$$\sum_{n=0}^{\infty} n \phi_n(t) z^{n-1} = \frac{1+a(p) k (t-1) z^{k-1}}{(1-z-a(p)(t-1)z^k)^2}\quad{\rm and }\quad \sum_{n=0}^{\infty}\phi_{n}^{\prime}(t) z^n = \frac{a(p) z^k}{(1-z-a(p)(t-1)z^k)^2}.$$
Multiplying the second expression with $(z+a(p)k(t-1)z^k)$ and comparing with the first expression, we get
\begin{align*}
(z+a(p)k(t-1)z^k) \sum_{n=0}^{\infty}\phi_{n}^{\prime}(t) z^n &=a(p) z^{k+1} \left(\frac{1+a(p)k(t-1)z^k}{(1-z-a(p)(t-1)z^k)^2}\right)= a(p) z^{k+1} \sum_{n=0}^{\infty} n \phi_n(t)z^{n-1}.
\end{align*}
Next, we compare the coefficients of $z^n$ and the required result follows.\\
(iii) Subtracting (ii) from (i) gives the result.\\
(iv) Using induction of $s$, for $2 \le s \le n$, we first prove
\begin{equation}
\phi_{n-s}^\prime(t) = - a(p) \sum_{u=0}^{s-2} \left(\frac{n-s+u+1}{k-1}\right) \left(\frac{k}{k-1}\right)^{u}\phi_{n-k+u-s+1}(t) + \left(\frac{k}{k-1}\right)^{s-1} \phi_{n-1}^\prime(t).\label{two:sss}
\end{equation}
Let $s=2$, replacing $n$ by $n-1$ in relation (iii), we have
\begin{align}
\phi_{n-2}^\prime(t) &= \phi_{n-1}^\prime(t) - a(p)\phi_{n-k-1}(t) - a(p)(t-1)\phi_{n-k-1}^{\prime}(t)\nonumber\\
&= \phi_{n-1}^\prime(t) - a(p)\phi_{n-k-1}(t) - \frac{1}{k-1}\left[(n-k) a(p) \phi_{n-k-1}(t)-\phi_{n-1}^\prime(t)\right]~~\{\rm by ~(i)\}\nonumber\\
&=- \left(\frac{n-1}{k-1}\right)a(p) \phi_{n-k-1}(t) + \frac{k}{k-1} \phi_{n-1}^{\prime}(t).\label{two:s=2}
\end{align}
Assume \eqref{two:sss} holds for $s=l$, then
$$\phi_{n-l}^\prime(t) = - a(p) \sum_{u=0}^{l-2} \left(\frac{n-l+u+1}{k-1}\right) \left(\frac{k}{k-1}\right)^{u}\phi_{n-k+u-l+1}(t) + \left(\frac{k}{k-1}\right)^{l-1} \phi_{n-1}^\prime(t).$$
For $s=l+1$, replacing $n$ by $n-1$ in above expression, we have
\begin{align*}
\phi_{n-(l+1)}^\prime(t) &= - a(p) \sum_{u=0}^{l-2} \left(\frac{n-l+u}{k-1}\right) \left(\frac{k}{k-1}\right)^{u}\phi_{n-k+u-l}(t) + \left(\frac{k}{k-1}\right)^{l-1} \phi_{n-2}^\prime(t)\\
&= - a(p) \sum_{u=0}^{l-2} \left(\frac{n-l+u}{k-1}\right) \left(\frac{k}{k-1}\right)^{u}\phi_{n-k+u-l}(t)\\
&~~~~~~~~~~~~~~~~~~~~~~~~~~~~~~~~~~~ + \left(\frac{k}{k-1}\right)^{l-1}\left[- \left(\frac{n-1}{k-1}\right)a(p) \phi_{n-k-1}(t) + \frac{k}{k-1} \phi_{n-1}^{\prime}(t)\right]\quad\{\rm by~\eqref{two:s=2}\}\\
&=- a(p) \sum_{u=0}^{l-1} \left(\frac{n-l+u}{k-1}\right) \left(\frac{k}{k-1}\right)^{u}\phi_{n-k+u-l}(t) + \left(\frac{k}{k-1}\right)^{l} \phi_{n-1}^\prime(t).
\end{align*}
Now, let $s=k$ in \eqref{two:sss} and replace right-hand side of relation (ii) for $\phi_{n-1}^{\prime}$. Then, the result follows by substituting $(n-k)$ with $n$.\qed\\
Next, we discuss the results related to the PMF of $B_{k_1,k_2}^n$. From \eqref{two:dpgf} (see Theorem $3.2$ of Dafnis {\em et al.} \cite{DAP} for details), the PMF satisfies the following relation
\begin{equation}
p_{m,n} = p_{m,n-1} + a(p) \left[p_{m-1,n-k}-p_{m,n-k}\right],\label{two:pmff}
\end{equation}
with initial conditions $p_{m,n}=\delta_{m,0}$ for $0 \le n < k$ and $p_{m,n}=0$ for $m <0$ or $m > \floor{n/k}$, where $\delta_{i,j}$ is the Kronecker delta function. Furthermore, \eqref{two:pmff} can be modified in the following form
\begin{equation}
p_{m,n} = p_{m,n-l} + a(p)\displaystyle{\sum_{s=0}^{l-1}\left[p_{m-1, n-k-s} - p_{m,n-k-s}\right]}, \quad \quad {\rm for}~ 1 \le l \le n-k+1.\label{two:relpmf}
\end{equation}
The following lemma gives an explicit form for PMF of $B_{k_1,k_2}^n$.

\begin{lemma}
The PMF of $B_{k_1,k_2}^n$ is given by
$$p_{m,n}={\mathbb P}\left(B_{k_1,k_2}^n=m \right)= \sum_{l=0}^{\floor{\frac{n-mk}{k}}} {n-(l+m)(k-1) \choose n-(l+m)k,~l,~m} (-1)^l a(p)^{l+m},\quad \quad m=0,1,\dotsc,\floor{n/k},$$
where $\displaystyle{{n \choose m_1,m_2,\dotsc,m_l}= \frac{n!}{{m_1}!~{m_2}!\dotsb{m_l}!}}$.
\end{lemma}
\noindent
{\bf Proof}. Multiplying both sides of \eqref{two:pmff} by $z_1^n z_2^m$, summing over $m$ and $n$ and rearranging the terms, we get
$$\left[1-z_1+a(p)z_1^k(1-z_2)\right]\sum_{n=0}^{\infty} \sum_{m=0}^{\infty} p_{m,n}z_1^n z_2^m = 1.$$
Therefore, for $|z_1-a(p)z_1^k(1-z_2)| < 1$, we have
$$\sum_{n=0}^{\infty} \sum_{m=0}^{\infty} p_{m,n}z_1^n z_2^m = \frac{1}{1\hspace{-0.1cm}-\hspace{-0.1cm}z_1+a(p)z_1^k(1\hspace{-0.1cm}-\hspace{-0.1cm}z_2)} = \sum_{n=0}^{\infty} z_1^n(1+a(p)z_1^{k-1}(z_2\hspace{-0.1cm}-\hspace{-0.1cm}1))^n = \sum_{n=0}^{\infty} \sum_{l=0}^{n} {n \choose l}a(p)^l z_1^{n+l(k-1)}(z_2\hspace{-0.1cm}-\hspace{-0.1cm}1)^l.$$
Interchanging the summation and adjusting the indices yield
$$\sum_{n=0}^{\infty} \sum_{m=0}^{\infty} p_{m,n}z_1^n z_2^m = \sum_{n=0}^{\infty} \sum_{l=0}^{\infty} {n+l \choose l}a(p)^l z_1^{n+lk} (z_2-1)^l=\sum_{n=0}^{\infty} \sum_{l=0}^{\infty} \sum_{m=0}^{l}{n+l \choose l}{l \choose m}(-1)^{l-m}a(p)^l z_1^{n+lk} z_2^m.$$
This leads to
$$\sum_{n=0}^{\infty} \sum_{m=0}^{\infty} p_{m,n}z_1^n z_2^m = \sum_{n=0}^{\infty} \sum_{l=0}^{\infty} \sum_{m=0}^{\infty}{n+l+m \choose n,~l,~m}(-1)^{l}a(p)^{l+m} z_1^{n+(l+m)k} z_2^m.$$
Replacing $n+(l+m)k$ by $n$, we get
$$\sum_{n=0}^{\infty} \sum_{m=0}^{\infty} p_{m,n}z_1^n z_2^m = \sum_{n=0}^{\infty}\sum_{m=0}^{\infty}\left( \sum_{l=0}^{\floor{\frac{n-mk}{k}}} {n-(l+m)(k-1)\choose n-(l+m)k,~l,~m}(-1)^{l}a(p)^{l+m} \right)z_1^{n} z_2^m.$$
Collecting the coefficients of $z_1^n z_2^m$ the result follows.\qed

\subsection{A Stein Operator for $B_{k_1,k_2}^n$}
The PGF approach is practically very useful for the distribution of runs as it is generally restricted to PGF. This approach is easy to apply whenever a relation between PGF and its derivative is known. So, to obtain Stein operator, we use PGF approach with the recursive relation (iv) derived in Lemma \ref{two:relpgf}. For more details and applications, we refer the reader to Upadhye {\em et al.} \cite{UCV} and Kumar and Upadhye \cite{KU}.

\begin{proposition}\label{two:ppp1}
A Stein operator for $B_{k_1,k_2}^n$ is given by
$${\cal A}_1g(m)=\left(\frac{n}{k}\hspace{-0.05cm}-\hspace{-0.05cm}m\right)\tilde{p}g(m+1)\hspace{-0.05cm}-\hspace{-0.05cm}(1-\tilde{p})mg(m)\hspace{-0.05cm}-\hspace{-0.05cm}a(p) \sum_{u=0}^{k-2} \left(\frac{n+u+1}{k-1}\right) \left(\frac{k}{k-1}\right)^{u}{\mathbb E}\big[g(B_{k_1,k_2}^{n\hspace{-0.01cm}-\hspace{-0.01cm}k\hspace{-0.01cm}+\hspace{-0.01cm}u\hspace{-0.01cm}+\hspace{-0.01cm}1}+1)|B_{k_1,k_2}^n=m\big].$$
\end{proposition}
\noindent
{\bf Proof}. From relation (iv) of Lemma \ref{two:relpgf}, we have
$$\displaystyle{\left[1+\tilde{p}(t-1)\right]\phi_{n}^{\prime}(t) =\frac{n}{k} \tilde{p} \phi_n(t) - a(p) \sum_{u=0}^{k-2} \left(\frac{n+u+1}{k-1}\right) \left(\frac{k}{k-1}\right)^{u}\phi_{n-k+u+1}(t)}.$$
This can be written as
\begin{align*}
\left[1+\tilde{p}(t-1)\right]&\sum_{m=0}^{\floor{n/k}}(m+1)p_{m+1,n}t^m=\frac{n}{k} \tilde{p}\sum_{m=0}^{\floor{n/k}}p_{m,n}t^m-a(p) \sum_{m=0}^{\floor{n/k}}\sum_{u=0}^{k-2} \left(\frac{n+u+1}{k-1}\right) \left(\frac{k}{k-1}\right)^{u}p_{m,n-k+u+1}t^m,
\end{align*}
where the summation of $m$ is independent of $u$ as $\floor{(n-k+u+1)/k}\le \floor{n/k}$ and $p_{m,n-k+u+1}$ is zero outside of $\floor{(n-k+u+1)/k}$. Collecting the coefficient of $t^m$, we get
\begin{align*}
\left(\frac{n}{k}-m\right)\tilde{p}~p_{m,n}-(1-\tilde{p})(m+1)p_{m+1,n}-a(p) \sum_{u=0}^{k-2} \left(\frac{n+u+1}{k-1}\right) \left(\frac{k}{k-1}\right)^{u}p_{m,n-k+u+1}=0.
\end{align*}
Let $g \in {\cal G}_{B_{k_1,k_2}^n}$, then
\begin{align*}
\sum_{m=0}^{\floor{n/k}}g(m+1)&\left[\left(\frac{n}{k}-m\right)\tilde{p}~p_{m,n}-(1-\tilde{p})(m+1)p_{m+1,n}-a(p) \sum_{u=0}^{k-2} \left(\frac{n+u+1}{k-1}\right) \left(\frac{k}{k-1}\right)^{u}p_{m,n-k+u+1}\right]=0.
\end{align*}
Rearranging the terms, we get
\begin{align*}
\sum_{m=0}^{\floor{n/k}}\left\{\left[\left(\frac{n}{k}-m\right)\tilde{p}g(m+1)\hspace{-0.07cm}-\hspace{-0.07cm}(1-\tilde{p})mg(m)\right]p_{m,n}\hspace{-0.07cm}-\hspace{-0.07cm}a(p) \sum_{u=0}^{k-2} \left(\frac{n+u+1}{k-1}\right) \left(\frac{k}{k-1}\right)^{u}\hspace{-0.15cm}g(m+1)p_{m,n-k+u+1}\right\}\hspace{-0.07cm}=\hspace{-0.07cm}0.
\end{align*}
Hence, a Stein operator for $B_{k_1,k_2}^n$ is given by
$${\cal A}_1g(m)=\left(\frac{n}{k}\hspace{-0.05cm}-\hspace{-0.05cm}m\right)\tilde{p}g(m+1)\hspace{-0.05cm}-\hspace{-0.05cm}(1-\tilde{p})mg(m)\hspace{-0.05cm}-a\hspace{-0.05cm}(p) \sum_{u=0}^{k-2} \left(\frac{n\hspace{-0.02cm}+\hspace{-0.02cm}u\hspace{-0.02cm}+\hspace{-0.02cm}1}{k-1}\right) \left(\frac{k}{k-1}\right)^{u}{\mathbb E}\big[g(B_{k_1,k_2}^{n-k+u+1}+1)|B_{k_1,k_2}^n=m\big].$$
This proves the result.\qed

\subsection{Pseudo-binomial Perturbation}
From Proposition \ref{two:ppp1}, it is clear that the Stein operator for $B_{k_1,k_2}^n$ can be seen as a perturbation of pseudo-binomial operator with parameter $n/k$ and $\tilde{p}$. Consider now pseudo-binomial distribution with general parameter $\alpha$ and $\check{p}$ as discussed in Section \ref{two:intro} and derive a Stein operator for $B_{k_1,k_2}^n$ by using recursive relations proved in Lemma \ref{two:relpgf} and parameter $\alpha$ and $\check{p}$. Later, using perturbation technique, we derive total variation distance between pseudo-binomial and $B_{k_1,k_2}^n$ by matching up to first two moments.\\
Recall that, from \eqref{two:MVAR}, the mean and variance of $B_{k_1,k_2}^n$ (with $p_l=p$, for all $l=1,2,\dotsc,n$) are as follows:
$${\mathbb E}\big(B_{k_1,k_2}^n\big)=(n-k+1)a(p)\quad {\rm and}\quad Var\big(B_{k_1,k_2}^n\big)=(n-k+1)(a(p)-(a(p))^2)-2\sum_{\substack{l<r\\r-l\le k-1}}a(p)^2.$$
Observe that
$$
2\hspace{-0.15cm}\sum_{\substack{l<r\\r-l\le k-1}}\hspace{-0.35cm}a(p)^2=2\hspace{-0.15cm}\sum_{l=1}^{n-k+1}\sum_{r=l+1}^{\min\{l+k-1,n-k+1\}}\hspace{-0.35cm}a(p)^2=\left\{
     \begin{array}{ll}
     0, & \text{if } n\le k\\
     (n-k)(n-k+1)a(p)^2, & \text{if } k < n \le 2k-1\\
     (2(n\hspace{-0.05cm}-\hspace{-0.05cm}2k\hspace{-0.05cm}+\hspace{-0.05cm}2)(k\hspace{-0.05cm}-\hspace{-0.05cm}1)+(k\hspace{-0.05cm}-\hspace{-0.05cm}1)(k\hspace{-0.05cm}-\hspace{-0.05cm}2))a(p)^2, & \text{if } n \ge 2k.
     \end{array}\right.
$$
Therefore, for $n\ge 2k$,
$${\mathbb E}\big(B_{k_1,k_2}^n\big)=(n-k+1)a(p)\quad {\rm and}\quad Var\big(B_{k_1,k_2}^n\big)=(n-k+1)a(p) +\left[(k-1)(3k-1)-(2k-1)n\right](a(p))^2.$$
For one-parameter approximation, we match mean of $B_{k_1,k_2}^n$ and $Z$ (with PMF given in \eqref{two:ppmf}) as follows:
\begin{equation}
\alpha \check{p} = (n-k+1) a(p). \label{two:mean}
\end{equation}
Here, matching can be done in two ways:
\begin{itemize}
\item[(a)] Fix $\alpha>0$ of our choice and $\check{p}=((n-k+1)/\alpha)a(p)$.
\item[(b)] Fix $0<\check{p}<1$ and $\alpha=((n-k+1)/\check{p})a(p)$.
\end{itemize}
For two-parameter approximation, we match first two moments as follows:
$$\alpha \check{p} = (n-k+1)a(p) \quad {\rm and} \quad \alpha \check{p} \check{q} = (n-k+1)a(p) +\left[(k-1)(3k-1)-(2k-1)n\right](a(p))^2.$$
This gives
\begin{equation}
\check{p} = \frac{\left[(2k-1)n-(k-1)(3k-1)\right]a(p)}{(n-k+1)}\quad {\rm and} \quad \alpha = \frac{(n-k+1)^2}{(2k-1)n-(k-1)(3k-1)}.\label{two:variance}
\end{equation}

\begin{theorem}\label{two:th1}
Let $n \ge 2 k \ge 4$ and $\tilde{p} < 1/2$ with \eqref{two:mean}, then
\begin{align*}
d_{TV}\left(B_{k_1,k_2}^n,~Z\right) &\le \frac{a(p)}{\floor{\alpha} \check{p}\check{q}}\left\{ \left( n(2k^{\star}\hspace{-0.1cm}-\hspace{-0.1cm}1)\hspace{-0.1cm}+\hspace{-0.1cm}k\hspace{-0.1cm}-\hspace{-0.1cm}1\right)\frac{|\tilde{p}-\check{p}|}{(1-2 \tilde{p})}\hspace{-0.1cm}+\hspace{-0.1cm}\left(n(k(k^{\star}\hspace{-0.1cm}-\hspace{-0.1cm}2)\hspace{-0.1cm}+\hspace{-0.1cm}1)-k(k\hspace{-0.1cm}-\hspace{-0.1cm}1)k^{\star}\hspace{-0.1cm}+\hspace{-0.1cm}3 k^2 \hspace{-0.1cm}-\hspace{-0.1cm}4 k \hspace{-0.1cm}+\hspace{-0.1cm}1\right)a(p)\right\},
\end{align*}
where $\displaystyle{k^{\star}=\left(\frac{k}{k-1}\right)^{k-1}}$.
\end{theorem}

\noindent
{\bf Proof}. First, we derive Stein operator for $B_{k_1,k_2}^n$ using PGF approach (Upadhye {\em et al.} \cite{UCV}). We have already shown in Lemma \ref{two:relpgf} that
$$\left[1+\tilde{p}(t-1)\right]\phi_{n}^{\prime}(t) =\frac{n}{k} \tilde{p} \phi_n(t) - a(p) \sum_{u=0}^{k-2} \left(\frac{n+u+1}{k-1}\right) \left(\frac{k}{k-1}\right)^{u}\phi_{n-k+u+1}(t).$$
Let $|t|< (1-\tilde{p})/\tilde{p}$, then this can be written as
\begin{align*}
\phi_n^\prime(t) &= \frac{1}{(1-\tilde{p} + \tilde{p} t)} \left(\frac{n}{k} \tilde{p} \phi_n(t) - a(p) \sum_{u=0}^{k-2} \left(\frac{n+u+1}{k-1}\right) \left(\frac{k}{k-1}\right)^{u}\phi_{n-k+u+1}(t)\right)\\
&= \frac{1}{1-\tilde{p}} \left(\sum_{m=0}^{\infty}(-1)^m \left(\frac{\tilde{p}}{1-\tilde{p}}\right)^m t^m\right) \left[\frac{n}{k} \tilde{p} \left(\sum_{m=0}^{\infty} p_{m,n} {\bf 1}\left(m \le \floor{\frac{n}{k}}\right)t^m\right) \right.\\
&~~~\left.- a(p) \sum_{u=0}^{k-2} \left(\frac{n+u+1}{k-1}\right) \left(\frac{k}{k-1}\right)^{u} \left(\sum_{m=0}^{\infty} p_{m,n-k+u+1} {\bf 1}\left(m \le \floor{\frac{n-k+u+1}{k}}\right)t^m\right)\right]\\
&= \frac{n \tilde{p}}{k (1-\tilde{p})} \sum_{m=0}^{\infty}\left( \sum_{l=0}^{m} p_{l,n} {\bf 1}\left(l \le \floor{\frac{n}{k}}\right)(-1)^{m-l} \left(\frac{\tilde{p}}{1-\tilde{p}}\right)^{m-l}\right)t^m- \frac{a(p)}{(1-\tilde{p})} \sum_{u=0}^{k-2} \left(\frac{n+u+1}{k-1}\right) \times\\
&~~~\left(\frac{k}{k-1}\right)^{u}\sum_{m=0}^{\infty}\left(\sum_{l=0}^{m}p_{l,n-k+u+1} {\bf 1}\left(l \le \floor{\frac{n-k+u+1}{k}}\right) (-1)^{m-l} \left(\frac{\tilde{p}}{1-\tilde{p}}\right)^{m-l}\right)t^m.
\end{align*}
Multiplying by $(\check{q}+\check{p}t)$ and comparing the coefficient of $t^m$, we get
\begin{align*}
&\check{q} (m+1) p_{m+1,n}~{\bf 1}\left(m \le \floor{\frac{n}{k}}-1\right) + \check{p} m p_{m,n} ~{\bf 1}\left(m \le \floor{\frac{n}{k}}\right)\\
& =  \frac{n \tilde{p}}{k (1-\tilde{p})} \left(\check{q} \sum_{l=0}^{m}p_{l,n}~ {\bf 1}\left(l \le \floor{\frac{n}{k}}\right) (-1)^{m-l} \left(\frac{\tilde{p}}{1-\tilde{p}}\right)^{m-l}\hspace{-0.2cm} + \check{p} \sum_{l=0}^{m-1}p_{l,n} ~{\bf 1}\left(l \le \floor{\frac{n}{k}}\right) (-1)^{m-l-1} \left(\frac{\tilde{p}}{1-\tilde{p}}\right)^{m-l-1}\right)\\
&~~-\frac{a(p)}{(1-\tilde{p})}\sum_{u=0}^{k-2} \left(\frac{n+u+1}{k-1}\right) \left(\frac{k}{k-1}\right)^{u}\left(\check{q}\sum_{l=0}^{m}p_{l,n-k+u+1} ~{\bf 1}\left(l \le \floor{\frac{n-k+u+1}{k}}\right) (-1)^{m-l} \left(\frac{\tilde{p}}{1-\tilde{p}}\right)^{m-l} \right.\\
&~~~\left. + \check{p} \sum_{l=0}^{m-1}p_{l,n-k+u+1}~ {\bf 1}\left(l \le \floor{\frac{n-k+u+1}{k}}\right) (-1)^{m-l-1} \left(\frac{\tilde{p}}{1-\tilde{p}}\right)^{m-l-1} \right),~~~~~~~~~~~~~~~~~~~~~~~~~~~~~~~~~~~~~~~~~~~~~~~~~~~~~~
\end{align*}
where $\check{q}=1-\check{p}$ as defined in \eqref{two:mean}. Rearranging the terms, we have
\begin{align*}
& \frac{n (\tilde{p}-\check{p})}{k (1-\tilde{p})}  \sum_{l=0}^{m}p_{l,n} ~{\bf 1}\left(l \le \floor{\frac{n}{k}}\right) (-1)^{m-l}\left(\frac{\tilde{p}}{1-\tilde{p}}\right)^{m-l}- \check{q} (m+1) p_{m+1,n} ~{\bf 1}\left(m \le \floor{\frac{n}{k}}-1\right)\\
&+\check{p}\left(\frac{n}{k}-m\right) p_{m,n}~ {\bf 1}\left(m \le \floor{\frac{n}{k}}\right)  - \frac{a(p)\check{p}}{\tilde{p}}\sum_{u=0}^{k-2} \left(\frac{n+u+1}{k-1}\right) \left(\frac{k}{k-1}\right)^{u}p_{m,n-k+u+1}~{\bf 1}\left(m \le \floor{\frac{n-k+u+1}{k}}\right) \\
&-\frac{a(p)(\tilde{p}-\check{p})}{\tilde{p}(1-\tilde{p})}\sum_{u=0}^{k-2} \left(\frac{n+u+1}{k-1}\right) \left(\frac{k}{k-1}\right)^{u}\sum_{l=0}^{m}p_{l,n-k+u+1} ~{\bf 1}\left(l \le \floor{\frac{n-k+u+1}{k}}\right) (-1)^{m-l} \left(\frac{\tilde{p}}{1-\tilde{p}}\right)^{m-l}= 0.
\end{align*}
Let $g \in {\cal G}_{B_{k_1,k_2}^n}$ as defined in \eqref{two:gx}, then
\begin{align*}
&\sum_{m=0}^{\infty}g(m+1)\left[ \frac{n (\tilde{p}-\check{p})}{k (1-\tilde{p})}  \sum_{l=0}^{m}p_{l,n} ~{\bf 1}\left(l \le \floor{\frac{n}{k}}\right) (-1)^{m-l}\left(\frac{\tilde{p}}{1-\tilde{p}}\right)^{m-l}- \check{q} (m+1) p_{m+1,n} ~{\bf 1}\left(m \le \floor{\frac{n}{k}}-1\right)\right.\\
&+\check{p}\left(\frac{n}{k}-m\right) p_{m,n}~ {\bf 1}\left(m \le \floor{\frac{n}{k}}\right)  - \frac{a(p)\check{p}}{\tilde{p}}\sum_{u=0}^{k-2} \left(\frac{n+u+1}{k-1}\right) \left(\frac{k}{k-1}\right)^{u}p_{m,n-k+u+1}~{\bf 1}\left(m \le \floor{\frac{n-k+u+1}{k}}\right) \\
&\left.-\frac{a(p)(\tilde{p}-\check{p})}{\tilde{p}(1-\tilde{p})}\sum_{u=0}^{k-2} \left(\frac{n+u+1}{k-1}\right) \left(\frac{k}{k-1}\right)^{u}\sum_{l=0}^{m}p_{l,n-k+u+1} ~{\bf 1}\left(l \le \floor{\frac{n-k+u+1}{k}}\right) (-1)^{m-l} \left(\frac{\tilde{p}}{1-\tilde{p}}\right)^{m-l}\right] \hspace{-0.1cm}= 0.
\end{align*}
Next, interchanging the sums lead to
\begin{align*}
&\sum_{m=0}^{\floor{n/k}}\left[ \check{p}\left(\frac{n}{k}-m\right) g(m+1) - \check{q} m g(m) \right] p_{m,n} + \frac{n (\tilde{p}-\check{p})}{k (1-\tilde{p})}  \sum_{l=0}^{\floor{n/k}}p_{l,n}\sum_{m=l}^{\infty} g(m+1)(-1)^{m-l}\left(\frac{\tilde{p}}{1-\tilde{p}}\right)^{m-l}\\
&-\frac{a(p)(\tilde{p}-\check{p})}{\tilde{p}(1-\tilde{p})}\sum_{u=0}^{k-2} \left(\frac{n+u+1}{k-1}\right) \left(\frac{k}{k-1}\right)^{u}\sum_{l=0}^{\floor{(n-k+u+1)/k}}p_{l,n-k+u+1}\sum_{m=l}^{\infty}g(m+1) (-1)^{m-l} \left(\frac{\tilde{p}}{1-\tilde{p}}\right)^{m-l}\\
&-a(p) \frac{\check{p}}{\tilde{p}}\sum_{u=0}^{k-2} \left(\frac{n+u+1}{k-1}\right) \left(\frac{k}{k-1}\right)^{u} \sum_{m=0}^{\floor{(n-k+u+1)/k}} g(m+1)p_{m,n-k+u+1} = 0.
\end{align*}
Substituting $m-l$ with $l$ and interchanging $m$ and $l$ (for second and third terms), this can be rewritten as
\begin{align*}
&\sum_{m=0}^{\floor{n/k}}\left[ \check{p}\left(\frac{n}{k}-m\right) g(m+1) - \check{q} m g(m) \right] p_{m,n} + \frac{n (\tilde{p}-\check{p})}{k (1-\tilde{p})}  \sum_{m=0}^{\floor{n/k}}p_{m,n}\sum_{l=0}^{\infty} g(m+l+1)(-1)^{l}\left(\frac{\tilde{p}}{1-\tilde{p}}\right)^{l}\\
&-\frac{a(p)(\tilde{p}-\check{p})}{\tilde{p}(1-\tilde{p})}\sum_{u=0}^{k-2} \left(\frac{n+u+1}{k-1}\right) \left(\frac{k}{k-1}\right)^{u}\sum_{m=0}^{\floor{(n-k+u+1)/k}}p_{m,n-k+u+1}\sum_{l=0}^{\infty}g(m+l+1) (-1)^{l} \left(\frac{\tilde{p}}{1-\tilde{p}}\right)^{l}\\
&-a(p) \frac{\check{p}}{\tilde{p}}\sum_{u=0}^{k-2} \left(\frac{n+u+1}{k-1}\right) \left(\frac{k}{k-1}\right)^{u} \sum_{m=0}^{\floor{(n-k+u+1)/k}} g(m+1)p_{m,n-k+u+1} = 0.
\end{align*}
Hence, Stein operator of $B_{k_1,k_2}^n$ is given by
\begin{align}
{\cal A}g(m) &=  \check{p}\left(\frac{n}{k}-m\right) g(m+1) - \check{q} m g(m) + \frac{n (\tilde{p}-\check{p})}{k (1-\tilde{p})}\sum_{l=0}^{\infty} g(m+l+1)(-1)^{l}\left(\frac{\tilde{p}}{1-\tilde{p}}\right)^{l}\nonumber\\
&-\frac{a(p)(\tilde{p}-\check{p})}{\tilde{p}(1-\tilde{p})}\sum_{u=0}^{k-2} \left(\frac{n+u+1}{k-1}\right) \left(\frac{k}{k-1}\right)^{u}{\mathbb E}\left(\sum_{l=0}^{\infty}g\left(B_{k_1,k_2}^{n-k+u+1}+l+1\right) (-1)^{l} \left(\frac{\tilde{p}}{1-\tilde{p}}\right)^{l}\Bigr|B_{k_1,k_2}^{n}=m\right)\nonumber\\
&-a(p) \frac{\check{p}}{\tilde{p}}\sum_{u=0}^{k-2} \left(\frac{n+u+1}{k-1}\right) \left(\frac{k}{k-1}\right)^{u} {\mathbb E} \left(g\left(B_{k_1,k_2}^{n-k+u+1}+1\right)\Bigr| B_{k_1,k_2}^{n}=m\right).\label{two:steinoper}
\end{align}
Next, introduce a parameter $\alpha > 0$ in \eqref{two:steinoper} and rewrite Stein operator as
\begin{align*}
\hat{\cal A}g(m) &=  \check{p}\left(\alpha-m\right) g(m+1) - \check{q} m g(m)+\left(\frac{n}{k}-\alpha\right)\check{p}g(m+1) + \frac{n (\tilde{p}-\check{p})}{k (1-\tilde{p})}\sum_{l=0}^{\infty} g(m+l+1)(-1)^{l}\left(\frac{\tilde{p}}{1-\tilde{p}}\right)^{l}\\
&-\frac{a(p)(\tilde{p}-\check{p})}{\tilde{p}(1-\tilde{p})}\sum_{u=0}^{k-2} \left(\frac{n+u+1}{k-1}\right) \left(\frac{k}{k-1}\right)^{u}{\mathbb E}\left(\sum_{l=0}^{\infty}g\left(B_{k_1,k_2}^{n-k+u+1}+l+1\right) (-1)^{l} \left(\frac{\tilde{p}}{1-\tilde{p}}\right)^{l}\Bigr|B_{k_1,k_2}^{n}=m\right)\\
&-a(p) \frac{\check{p}}{\tilde{p}}\sum_{u=0}^{k-2} \left(\frac{n+u+1}{k-1}\right) \left(\frac{k}{k-1}\right)^{u} {\mathbb E} \left(g\left(B_{k_1,k_2}^{n-k+u+1}+1\right)\Bigr| B_{k_1,k_2}^{n}=m\right)\\
&=\hat{\cal A}_0 g(m)+ \hat{\cal U}g(m),
\end{align*}
where $\hat{\cal A}_0$ is Stein operator of pseudo-binomial with parameter $(\alpha, \check{p})$ with $\check{p}$ defined as in \eqref{two:mean} and $\hat{\cal U}$ is the perturbed operator. Now taking the expectation w.r.t. $B_{k_1,k_2}^n$, we get
\begin{align*}
{\mathbb E}\left[\hat{\cal U}g\left(B_{k_1,k_2}^n\right)\right] &=\sum_{m=0}^{\floor{n/k}}\left(\frac{n}{k}-\alpha\right)\check{p}g(m+1)p_{m,n} + \frac{n (\tilde{p}-\check{p})}{k (1-\tilde{p})}  \sum_{m=0}^{\floor{n/k}}p_{m,n}\sum_{l=0}^{\infty} g(m+l+1)(-1)^{l}\left(\frac{\tilde{p}}{1-\tilde{p}}\right)^{l}\\
&~~~-\frac{a(p)(\tilde{p}\hspace{-0.07cm}-\hspace{-0.07cm}\check{p})}{\tilde{p}(1-\tilde{p})}\sum_{u=0}^{k-2} \left(\frac{n\hspace{-0.07cm}+\hspace{-0.07cm}u\hspace{-0.07cm}+\hspace{-0.07cm}1}{k-1}\right) \left(\frac{k}{k-1}\right)^{u}\sum_{m=0}^{\floor{n-k+u+1/k}}\hspace{-0.6cm}p_{m,n-k+u+1}\sum_{l=0}^{\infty}g(m+l+1) (-1)^{l} \left(\frac{\tilde{p}}{1-\tilde{p}}\right)^{l}\nonumber\\
&~~~-a(p) \frac{\check{p}}{\tilde{p}}\sum_{u=0}^{k-2} \left(\frac{n+u+1}{k-1}\right) \left(\frac{k}{k-1}\right)^{u} \sum_{m=0}^{\floor{n-k+u+1/k}} g(m+1)p_{m,n-k+u+1}.
\end{align*}
It is known that
\begin{equation}
g(m+l+1) = \sum_{j=1}^{l}\Delta g(m+j) + g(m+1). \label{two:delta}
\end{equation}
Using (\ref{two:delta}) and replacing $\floor{(n-k+u+1)/k}$ by $\floor{n/k}$ as $p_{m,n-k+u+1}$ is zero outside of its range, the expression becomes
\begin{align}
{\mathbb E}\left[\hat{\cal U}g\left(B_{k_1,k_2}^n\right)\right] &=\frac{n }{k }(\tilde{p}-\check{p})  \sum_{m=0}^{\floor{n/k}}p_{m,n}\sum_{j=1}^{\infty} \Delta g(m+j)(-1)^{j}\left(\frac{\tilde{p}}{1-\tilde{p}}\right)^{j}+\left(\frac{n}{k}\tilde{p}-\alpha\check{p}\right)\sum_{m=0}^{\floor{n/k}} g(m+1) p_{m,n}\nonumber\\
&~~~-\frac{a(p)}{\tilde{p}}(\tilde{p}-\check{p})\sum_{u=0}^{k-2} \left(\frac{n+u+1}{k-1}\right) \left(\frac{k}{k-1}\right)^{u}\sum_{m=0}^{\floor{n/k}}p_{m,n-k+u+1}\sum_{j=1}^{\infty}\Delta g(m+j)(-1)^{j} \left(\frac{\tilde{p}}{1-\tilde{p}}\right)^{j}\nonumber\\
&~~~-a(p)\sum_{u=0}^{k-2} \left(\frac{n+u+1}{k-1}\right) \left(\frac{k}{k-1}\right)^{u}\sum_{m=0}^{\floor{n/k}}p_{m,n-k+u+1}g(m+1).\label{two:exp0}
\end{align}
Observe that
\begin{align*}
\sum_{u=0}^{k-2} \left(\frac{n+u+1}{k-1}\right) \left(\frac{k}{k-1}\right)^{u}=n\left(\frac{k}{k-1}\right)^{k-1}-n+k-1=\frac{n}{k}\tilde{p}-\alpha\check{p}.
\end{align*}
Therefore, the last two terms of \eqref{two:exp0} can be written as
\begin{align}
\left(\frac{n}{k}\tilde{p}-\alpha\check{p}\right)\sum_{m=0}^{\floor{n/k}} g(m+1) p_{m,n}&-a(p)\sum_{u=0}^{k-2} \left(\frac{n+u+1}{k-1}\right) \left(\frac{k}{k-1}\right)^{u}\sum_{m=0}^{\floor{n/k}}p_{m,n-k+u+1}g(m+1)\nonumber\\
&= a(p)\sum_{u=0}^{k-2} \left(\frac{n+u+1}{k-1}\right) \left(\frac{k}{k-1}\right)^{u} \sum_{m=0}^{\floor{n/k}} \left[p_{m,n}-p_{m,n-k+u+1}\right]g(m+1).\nonumber
\end{align}
Now, using \eqref{two:relpmf}, for $l=k-u-1$, we get
\begin{align}
\left(\frac{n}{k}\tilde{p}-\alpha\check{p}\right)\sum_{m=0}^{\floor{n/k}} g(m+1) p_{m,n}&-a(p)\sum_{u=0}^{k-2} \left(\frac{n+u+1}{k-1}\right) \left(\frac{k}{k-1}\right)^{u}\sum_{m=0}^{\floor{n/k}}p_{m,n-k+u+1}g(m+1)\nonumber\\
&= a(p)\sum_{u=0}^{k-2} \left(\frac{n+u+1}{k-1}\right) \left(\frac{k}{k-1}\right)^{u} \sum_{m=0}^{\floor{n/k}} \sum_{s=0}^{k-u-2} \left[p_{m-1,n-k-s} - p_{m, n-k-s}\right]g(m+1)\nonumber\\
&= (a(p))^2\sum_{u=0}^{k-2} \left(\frac{n+u+1}{k-1}\right) \left(\frac{k}{k-1}\right)^{u} \sum_{m=0}^{\floor{n/k}} \sum_{s=0}^{k-u-2} \Delta g(m+1)p_{m,n-k-s}. \label{two:deltaeq4}
\end{align}
Next, combining (\ref{two:exp0}) and (\ref{two:deltaeq4}), we see that
\begin{align}
{\mathbb E} \left[\hat{\cal U}g\left(B_{k_1,k_2}^n\right)\right]&= \frac{n }{k } (\tilde{p}-\check{p}) \sum_{m=0}^{\floor{n/k}}p_{m,n} \sum_{j=1}^{\infty} \Delta g(m+j)(-1)^{j}\left(\frac{\tilde{p}}{1-\tilde{p}}\right)^{j}\nonumber\\
&~~-\frac{a(p)(\tilde{p}-\check{p})}{\tilde{p}}\sum_{u=0}^{k-2} \left(\frac{n+u+1}{k-1}\right) \left(\frac{k}{k-1}\right)^{u}\sum_{m=0}^{\floor{n/k}}p_{m,n-k+u+1}\sum_{j=1}^{\infty}\Delta g(m+j) (-1)^{j} \left(\frac{\tilde{p}}{1-\tilde{p}}\right)^{j}\nonumber\\
&~~+(a(p))^2\sum_{u=0}^{k-2} \left(\frac{n+u+1}{k-1}\right) \left(\frac{k}{k-1}\right)^{u} \sum_{m=0}^{\floor{n/k}} \sum_{s=0}^{k-u-2} \Delta g(m+1)p_{m,n-k-s}.\label{two:oneexp}
\end{align}
Therefore, for $\tilde{p} < 1/2$ and $g \in {\cal G}_{Z} \cap {\cal G}_{B_{k_1,k_2}^n}$, and using \eqref{two:bound}, we have
\begin{align*}
\left|{\mathbb E}\left[\hat{\cal U} g\left(B_{k_1,k_2}^n\right)\right]\right| &\le \|\Delta g\|\left[ \frac{n \tilde{p}}{k (1-2 \tilde{p})} |\tilde{p}-\check{p}|+\frac{a(p)}{(1-2 \tilde{p})}\left|\tilde{p}-\check{p}\right|\sum_{u=0}^{k-2} \left(\frac{n+u+1}{k-1}\right) \left(\frac{k}{k-1}\right)^{u}\right.\\
&~~~~~~~~~~~~~~~~~~~~~~~~~~~~~~~~~~~~~~~~~~~~~~~~~~~~~~~\left.+(a(p))^2\sum_{u=0}^{k-2} \left(\frac{n+u+1}{k-1}\right) \left(\frac{k}{k-1}\right)^{u}(k-u-1)\right]\\
&\le\frac{a(p)}{\floor{\alpha} \check{p}\check{q}}\left\{ \left( n(2k^{\star}\hspace{-0.1cm}-\hspace{-0.1cm}1)+k\hspace{-0.1cm}-\hspace{-0.1cm}1\right)\frac{|\tilde{p}-\check{p}|}{(1-2 \tilde{p})}+\left(n(k(k^{\star}\hspace{-0.1cm}-\hspace{-0.1cm}2)+1)\hspace{-0.1cm}-k(k\hspace{-0.1cm}-\hspace{-0.1cm}1)k^{\star}+3 k^2 \hspace{-0.1cm}-\hspace{-0.1cm}4 k +1\right)a(p)\right\}.
\end{align*}
This proves the result.\qed\\
\noindent
Next, for two-parameter approximation, we first generalize Theorem $4.7$ of Upadhye {\em et al.} \cite{UCV} for ($k_1,k_2$)-runs. The proof of Proposition \ref{two:peo} follows in a similar spirit of Theorem $4.7$ of Upadhye {\em et al.} \cite{UCV} and Lemma $2.1$ of Xia and Zhang \cite{WX}.

\begin{proposition}\label{two:peo}
Let $\{I_l\}$ and $B_{k_1,k_2}^n$ as defined in \eqref{two:ij} with $na(p) \ge 8$. Then
$$d_{TV}\left(B_{k_1,k_2}^n,~B_{k_1,k_2}^n+1\right) \le \min\left\{1,M(n)\right\},$$
where $M(n)=\frac{72(1-(2k-1)a(p))}{na(p)}+\sqrt{\frac{2}{\pi}} \left(\frac{1}{4}+n a(p)\left(1-a(p)\right)\right)^{-1/2}$.
\end{proposition}
\noindent
{\bf Proof}. Let us construct another version of $\{I_l\}$, denoted by $\{\tilde{I_l}\}$ such that $\tilde{B}_{k_1,k_2}^n=\sum_{l=1}^{n-k+1}\tilde{I_l}$ then it is enough to prove $d_{TV}\big(\tilde{B}_{k_1,k_2}^n,~B_{k_1,k_2}^n+1\big) \le \min\left\{1,M(n)\right\}$. Define stopping times
\vspace{-0.2cm}
$$\rho_m = \min\{l > \rho_{m-1} |I_l = 1\} \quad {\rm and} \quad \rho_0 = 0.$$
\vspace{-0.2cm}
Also, define $T_m = \rho_m - \rho_{m-1}$, then from Huang and Tsai \cite{HT} $T_m$’s are iid, with iid copy $T$, having PGF
$${\mathbb E}(z^T) = \frac{a(p)z^k}{1-z+a(p)z^k}.$$
It is clear that that ${\mathbb E}T = 1/a(p)$ and $Var(T) =\frac{1-(2k-1)a(p)}{(a(p))^2}$. Note that $\rho_m =\sum_{j=1}^{m}T_j$ is the waiting time for $m$th occurrence of $(k_1, k_2 )$-event. Then, the average number of occurrences of ($k_1,k_2$)-events is $\frac{n}{{\mathbb E}T} = na(p)$. Let $\upsilon = \floor{na(p)} +1$, then $\rho_\upsilon =\sum_{m=1}^{\upsilon}T_m$ and by Corollary $1.6$ of Mattner and Roos \cite{MR}, we have
$$d_{TV}\left(\rho_\upsilon,~\rho_{\upsilon}+1\right) \le \sqrt{\frac{2}{\pi}} \left(\frac{1}{4}+\sum_{m=1}^{\upsilon}\left(1-d_{TV}\left(T_m,T_m+1\right)\right)\right)^{-1/2}.$$
Now, it is known that $d_{T V} (T , T + 1) = a(p)$ which implies
$$d_{TV}\left(\rho_\upsilon,~\rho_{\upsilon}+1\right) \le \sqrt{\frac{2}{\pi}} \left(\frac{1}{4}+\upsilon\left(1-a(p)\right)\right)^{-1/2} \le \sqrt{\frac{2}{\pi}} \left(\frac{1}{4}+n a(p)\left(1-a(p)\right)\right)^{-1/2}.$$
Define maximal coupling similar to Xia and Zhang \cite{XZ} p.-1339 (see also, Barbour {\em et al.} \cite{BHJ} and Wang \cite{W}).
\begin{equation}
{\mathbb P}\left(\rho_\upsilon \neq \tilde{\rho}_\upsilon+1\right) = d_{TV}\left(\rho_\upsilon,~\rho_{\upsilon}+1\right) \le \sqrt{\frac{2}{\pi}} \left(\frac{1}{4}+n a(p)\left(1-a(p)\right)\right)^{-1/2}.\label{two:first}
\end{equation}
Now, let $\tilde{\rho}_\upsilon = \sum_{m=1}^{\upsilon} \tilde{T}_m$ such that $\tilde{T}_m$'s are iid and $\tilde{\rho}_m=\tilde{\rho}_{m-1} +\tilde{T}_m$ with $\tilde{\rho}_0=0$. Define now
\[
\tilde{I}_l =\left\{
     \begin{array}{ll}
     0, & \text{$\tilde{\rho}_{m-1} < l < \tilde{\rho}_{m} \quad 1 \le m \le \upsilon,$}\\
     1, & \text{$\tilde{\rho}_{m} =l \quad 1 \le m \le \upsilon,$}\\
     I_l, & \text{$\tilde{\rho}_{\upsilon} < l.$}
     \end{array}\right.
\]
Then, for $\rho_\upsilon \le n$ and ${\rho}_\upsilon = \tilde{\rho}_\upsilon + 1$, we have ${\tilde{B}_{k_1,k_2}^n}= B_{k_1,k_2}^n+ 1$. Hence,
\begin{equation}
{\mathbb P }\left({\tilde{B}_{k_1,k_2}^n} \neq {B_{k_1,k_2}^n+1} \right) \le {\mathbb P} (\rho_\upsilon > n) + {\mathbb P} (\rho_\upsilon \neq \tilde{\rho}_\upsilon+1).\label{two:second}
\end{equation}
From Chebyshev's inequality, we have
$${\mathbb P}(\rho_\upsilon>n) \le \frac{Var(\rho_\upsilon)}{(n - {\mathbb E}\rho_\upsilon)^2}.$$
We have already seen that
$${\mathbb E}\rho_\upsilon = \frac{\upsilon}{a(p)} \quad {\rm and}\quad Var(\rho_\upsilon) = \frac{\upsilon(1+(1-2k)a(p))}{(a(p))^2}.$$
Without loss of generality, let $na(p) \ge 8$, then
\begin{align}
{\mathbb P}(\rho_\upsilon>n) &\le \frac{\upsilon (a(p)(1-2k)+1)}{(na(p)-\upsilon)^2}\nonumber\\
&\le \frac{1.125 (a(p)(1-2k)+1)}{na(p)(0.125)^2}\nonumber\\
&= \frac{72 (a(p)(1-2k)+1)}{na(p)}= \frac{M}{n}.\label{two:third}
\end{align}
Combining \eqref{two:first}, \eqref{two:second} and \eqref{two:third} with $d_{TV}\big(\tilde{B}_{k_1,k_2}^n,~B_{k_1,k_2}^n+1\big) \le1$, we get the required result.\qed

\noindent
Next, we derive the total variation distance between pseudo-binomial and $B_{k_1,k_2}^n$ by matching first two moments as discussed in \eqref{two:variance}.
\begin{theorem}\label{two:th2}
Let $n \ge 3 k \ge 6$, $na(p) \ge 8$ and $\tilde{p} < 1/2$ with \eqref{two:variance}, then
\begin{align*}
d_{TV}\left(B_{k_1,k_2}^n,~Z\right) & \le \frac{2(a(p))^2}{\floor{\alpha}\check{p}\check{q}}\left\{|\tilde{p}\hspace{-0.07cm}-\hspace{-0.07cm}\check{p}|\left(\frac{k k^{\star}\left(n(2 k^{\star}\hspace{-0.1cm}-\hspace{-0.1cm}1)\hspace{-0.07cm}+\hspace{-0.07cm}k\hspace{-0.07cm}-\hspace{-0.07cm}1\right)}{1\hspace{-0.07cm}-\hspace{-0.07cm}2\tilde{p}}\hspace{-0.07cm}+\hspace{-0.07cm}\left(n(k(k^{\star}\hspace{-0.1cm}-\hspace{-0.1cm}2)+1)-k(k-1)k^{\star}+3 k^2-4k+1\right)\right)\right.\\
&~~~+\left(n((2k-1)k k^{\star}-\frac{9}{2}k(k-1)-1)-kk^* (3k-1)(k-1)+\left.\frac{17k^3-30k^2+15k-2}{2}\right)a(p)\right\}\\
&~~~~~~~~\times\left\{1 \wedge \left(\frac{M}{n-3k+3}+\sqrt{\frac{2}{\pi}} \left(\frac{1}{4}+(n-3k+3) a(p)\left(1-a(p)\right)\right)^{-1/2}\right)\right\}.
\end{align*}
\end{theorem}
\noindent
{\bf Proof}. From \eqref{two:oneexp}, we have
\begin{align}
{\mathbb E} \left[\hat{\cal U}g\left(B_{k_1,k_2}^n\right)\right]&= \frac{n }{k } (\tilde{p}-\check{p}) \sum_{m=0}^{\floor{n/k}}p_{m,n} \sum_{j=1}^{\infty} \Delta g(m+j)(-1)^{j}\left(\frac{\tilde{p}}{1-\tilde{p}}\right)^{j}\nonumber\\
&~~-a(p)\left(1-\frac{\check{p}}{\tilde{p}}\right)\sum_{u=0}^{k-2} \left(\frac{n+u+1}{k-1}\right) \left(\frac{k}{k-1}\right)^{u}\sum_{m=0}^{\floor{n/k}}p_{m,n-k+u+1}\sum_{j=1}^{\infty}\Delta g(m+j) (-1)^{j} \left(\frac{\tilde{p}}{1-\tilde{p}}\right)^{j}\nonumber\\
&~~+(a(p))^2\sum_{u=0}^{k-2} \left(\frac{n+u+1}{k-1}\right) \left(\frac{k}{k-1}\right)^{u} \sum_{m=0}^{\floor{n/k}} \sum_{s=0}^{k-u-2} \Delta g(m+1)p_{m,n-k-s}.\label{two:oneexp2}
\end{align}
Again, from Newton expansion (see Barbour {\em et al.} \cite{BCX}), we have
\begin{equation}
\Delta g(m+j) = \sum_{l=1}^{j-1} \Delta^2 g(m+l) + \Delta g(m+1).\label{two:ddelta}
\end{equation}
Using \eqref{two:ddelta} in \eqref{two:oneexp2}, we get
\begin{align}
{\mathbb E}\left[\hat{\cal U}g\left(B_{k_1,k_2}^n\right)\right] &=-\frac{n }{k }\tilde{p}(\tilde{p}-\check{p})  \sum_{m=0}^{\floor{n/k}}p_{m,n}\sum_{l=1}^{\infty} \Delta^2 g(m+l) (-1)^{l}\left(\frac{\tilde{p}}{1-\tilde{p}}\right)^{l}-\frac{n }{k }\tilde{p}(\tilde{p}-\check{p})  \sum_{m=0}^{\floor{n/k}}p_{m,n}\Delta g(m+1)\nonumber\\
&~~~+a(p)(\tilde{p}-\check{p})\sum_{u=0}^{k-2} \left(\frac{n+u+1}{k-1}\right) \left(\frac{k}{k-1}\right)^{u}\sum_{m=0}^{\floor{n/k}}p_{m,n-k+u+1}\sum_{l=1}^{\infty} \Delta^2 g(m+l) (-1)^{l} \left(\frac{\tilde{p}}{1-\tilde{p}}\right)^{l}\nonumber\\
&~~~+a(p)(\tilde{p}-\check{p})\sum_{u=0}^{k-2} \left(\frac{n+u+1}{k-1}\right) \left(\frac{k}{k-1}\right)^{u}\sum_{m=0}^{\floor{n/k}}p_{m,n-k+u+1}\Delta g(m+1)\nonumber\\
&~~~+a(p)^2\sum_{u=0}^{k-2}\left(\frac{n+u+1}{k-1}\right)\left(\frac{k}{k-1}\right)^{u}\sum_{m=0}^{\floor{n/k}}\sum_{s=0}^{k-u-2}\Delta g(m+1)p_{m,n-k-s}. \label{two:exp2}
\end{align}
Observe that
\begin{align*}
a(p)(\tilde{p}-\check{p})\sum_{u=0}^{k-2} \left(\frac{n+u+1}{k-1}\right) \left(\frac{k}{k-1}\right)^{u}+a(p)^2\sum_{u=0}^{k-2}(k-u-1)\left(\frac{n+u+1}{k-1}\right)\left(\frac{k}{k-1}\right)^{u}=\frac{n }{k }\tilde{p}(\tilde{p}-\check{p}).
\end{align*}
Combining the terms involving $\Delta g(m+1)$ in \eqref{two:exp2}, we get
\begin{align*}
{\mathbb E}\left[\hat{\cal U}g\left(B_{k_1,k_2}^n\right)\right] &=-\frac{n }{k }\tilde{p}(\tilde{p}-\check{p})  \sum_{m=0}^{\floor{n/k}}p_{m,n}\sum_{l=1}^{\infty} \Delta^2 g(m+l) (-1)^{l}\left(\frac{\tilde{p}}{1-\tilde{p}}\right)^{l}\nonumber\\
&~~~+a(p)(\tilde{p}-\check{p})\sum_{u=0}^{k-2} \left(\frac{n+u+1}{k-1}\right) \left(\frac{k}{k-1}\right)^{u}\sum_{m=0}^{\floor{n/k}}p_{m,n-k+u+1}\sum_{l=1}^{\infty} \Delta^2 g(m+l) (-1)^{l} \left(\frac{\tilde{p}}{1-\tilde{p}}\right)^{l}\nonumber\\
&~~~-a(p)(\tilde{p}-\check{p})\sum_{u=0}^{k-2} \left(\frac{n+u+1}{k-1}\right) \left(\frac{k}{k-1}\right)^{u}\sum_{m=0}^{\floor{n/k}}\Delta g(m+1)[p_{m,n}-p_{m,n-k+u+1}]
\end{align*}
\begin{align*}
&~~~~~~-a(p)^2\sum_{u=0}^{k-2}\left(\frac{n+u+1}{k-1}\right)\left(\frac{k}{k-1}\right)^{u}\sum_{m=0}^{\floor{n/k}}\sum_{s=0}^{k-u-2}\Delta g(m+1)[p_{m,n}-p_{m,n-k-s}].
\end{align*}
Using \eqref{two:relpmf} for the last two terms, we have
\begin{align*}
{\mathbb E}\left[\hat{\cal U}g\left(B_{k_1,k_2}^n\right)\right] &=-\frac{n }{k }\tilde{p}(\tilde{p}-\check{p})  \sum_{m=0}^{\floor{n/k}}p_{m,n}\sum_{l=1}^{\infty} \Delta^2 g(m+l) (-1)^{l}\left(\frac{\tilde{p}}{1-\tilde{p}}\right)^{l}\\
&~~~+a(p)(\tilde{p}-\check{p})\sum_{u=0}^{k-2} \left(\frac{n+u+1}{k-1}\right) \left(\frac{k}{k-1}\right)^{u}\sum_{m=0}^{\floor{n/k}}p_{m,n-k+u+1}\sum_{l=1}^{\infty} \Delta^2 g(m+l) (-1)^{l} \left(\frac{\tilde{p}}{1-\tilde{p}}\right)^{l}\\
&~~~-(a(p))^3\sum_{u=0}^{k-2} \left(\frac{n\hspace{-0.07cm}+\hspace{-0.07cm}u\hspace{-0.07cm}+\hspace{-0.07cm}1}{k-1}\right) \left(\frac{k}{k-1}\right)^{u}\sum_{m=0}^{\floor{n/k}} \sum_{s=0}^{k-u-2} \sum_{t=0}^{k+s-1} \Delta g(m+1) [p_{m-1,n-k-t}-p_{m,n-k-t}]\\
&~~~-(a(p))^2(\tilde{p}\hspace{-0.07cm}-\hspace{-0.07cm}\check{p})\sum_{u=0}^{k-2} \left(\frac{n\hspace{-0.07cm}+\hspace{-0.07cm}u\hspace{-0.07cm}+\hspace{-0.07cm}1}{k-1}\right) \left(\frac{k}{k-1}\right)^{u}\sum_{m=0}^{\floor{n/k}}\sum_{s=0}^{k-u-2}\Delta g(m+1)[p_{m-1,n-k-s}\hspace{-0.07cm}-\hspace{-0.07cm}p_{m,n-k-s}].
\end{align*}
Hence, for $g \in {\cal G}_Z \cap {\cal G}_{B_{k_1,k_2}^n}$, and using \eqref{two:bound} and Lemma \ref{two:peo}, we get
\vspace{-0.25cm}
\begin{align*}
\left|{\mathbb E}\left[\hat{\cal U}g\left(B_{k_1,k_2}^n\right)\right]\right| &\le 2 \|\Delta g\|\left\{\frac{n \tilde{p}^2|\tilde{p}-\check{p}|}{k (1-2\tilde{p}) } d_{TV}\left(B_{k_1,k_2}^n,~B_{k_1,k_2}^n+1\right) \right.\\
&~~~+\frac{a(p)\tilde{p}|\tilde{p}-\check{p}|}{1-2\tilde{p}}\sum_{u=0}^{k-2} \left(\frac{n+u+1}{k-1}\right) \left(\frac{k}{k-1}\right)^{u} d_{TV}\left(B_{k_1,k_2}^{n-k+u+1},B_{k_1,k_2}^{n-k+u+1}+1\right)\\
&~~~+(a(p))^3\sum_{u=0}^{k-2} \left(\frac{n+u+1}{k-1}\right) \left(\frac{k}{k-1}\right)^{u} \sum_{s=0}^{k-u-2} \sum_{t=0}^{k+s-1} d_{TV}\left(B_{k_1,k_2}^{n-k-t},~B_{k_1,k_2}^{n-k-t}+1\right)\\
&~~~\left.+(a(p))^2|\tilde{p}-\check{p}|\sum_{u=0}^{k-2} \left(\frac{n+u+1}{k-1}\right) \left(\frac{k}{k-1}\right)^{u}\sum_{s=0}^{k-u-2}d_{TV}\left(B_{k_1,k_2}^{n-k-s},~B_{k_1,k_2}^{n-k-s}+1\right)\right\}\\
& \le \frac{2(a(p))^2}{\floor{\alpha}\check{p}\check{q}}\left\{|\tilde{p}\hspace{-0.07cm}-\hspace{-0.07cm}\check{p}|\left(\frac{k k^{\star}\left(n(2 k^{\star}\hspace{-0.1cm}-\hspace{-0.1cm}1)\hspace{-0.07cm}+\hspace{-0.07cm}k\hspace{-0.07cm}-\hspace{-0.07cm}1\right)}{1\hspace{-0.07cm}-\hspace{-0.07cm}2\tilde{p}}\hspace{-0.07cm}+\hspace{-0.07cm}\left(n(k(k^{\star}\hspace{-0.1cm}-\hspace{-0.1cm}2)+1)-k(k-1)k^{\star}+3 k^2-4k+1\right)\right)\right.\\
&~~~+\left(n((2k-1)k k^{\star}-\frac{9}{2}k(k-1)-1)-kk^* (3k-1)(k-1)+\left.\frac{17k^3-30k^2+15k-2}{2}\right)a(p)\right\}\\
&~~~~~~~~\times\left\{ 1 \wedge \left(\frac{M}{n-3k+3}+\sqrt{\frac{2}{\pi}} \left(\frac{1}{4}+(n-3k+3) a(p)\left(1-a(p)\right)\right)^{-1/2}\right)\right\}.
\end{align*}
This proves the result.\qed

\section{Approximation Results for Independent Trials}\label{two:ARit}
In this section, we obtain the approximation results in total variation distance for $(k_1,k_2)$-runs arising from a sequence of independent Bernoulli trials by matching up to two moments.\\
As discussed in Section \ref{two:intro}, the distribution of $B_{k_1,k_2}^n$ becomes intractable whenever trials are non-identical. Also, the PGF of $B_{k_1,k_2}^n$ for independent Bernoulli trials (which can be derived (see Fu and Koutras \cite{FK})) can not be written in compact form. Hence, it is difficult to apply PGF approach to obtain a Stein operator. However, we can modify and use the technique discussed by Wang and Xia \cite{WX} for $(k_1,k_2)$-runs, by considering the Stein operator of pseudo-binomial distribution and obtain bounds in total variation distance. \\
Recall that $\xi_1,\xi_2,\dotsc,\xi_n$ are independent Bernoulli trials with probability ${\mathbb P}(\xi_l=1)=p_l=1-{\mathbb P}(\xi_l=0)$, $l=1,2,\dotsc,n$ and $B_{k_1,k_2}^n=\sum_{l=1}^{n-k_1-k_2+1} I_l$, where $I_l=(1-\xi_l)\dotsb(1-\xi_{l+k_1-1})\xi_{l+k_1}\dotsb\xi_{l+k_1+k_2-1}$. Define
$$W_l=B_{k_1,k_2}^n-I_l, \quad\quad X_l = B_{k_1,k_2}^n-\sum_{|r-l|\le k_1+k_2-1}I_r,\quad\quad D_{l,u}=\sum_{r=l-k_1-k_2+1}^{u}I_r\quad{\rm and}\quad \tilde{D}_{l,u}=\sum_{\substack{r=l-k_1-k_2+1\\r\neq l}}^{u}I_r .$$
Observe that $X_l$ and $I_l$ are independent. For one-parameter approximation, we match the first moment of $B_{k_1,k_2}^n$ and $Z$ as follows:
\begin{equation}
\alpha \check{p}=\sum_{l=1}^{n-k_1-k_2+1}{\mathbb E}(I_l)=\sum_{l=1}^{n-k_1-k_2+1}a(p_l).\label{two:mean2}
\end{equation}

\begin{theorem}\label{two:Ith1}
Let $n\ge 2(k_1+k_2)$ with \eqref{two:mean2}, then
$$d_{TV}\left(B_{k_1,k_2}^n,Z\right) \le \frac{1}{\floor{\alpha} \check{p}\check{q}}\sum_{l=1}^{n-k_1-k_2+1}a(p_l)\left\{\sum_{|u-l|\le k_1+k_2-1} a(p_u)+\check{p}\right\}.$$
\end{theorem}
\noindent
{\bf Proof}. Recall that, from \eqref{two:sboper}, a Stein operator for $Z$ is given by
$${\cal A}_0 g(m)=(\alpha-m)\check{p}g(m+1)-m\check{q}g(m).$$
Replacing $m$ with $B_{k_1,k_2}^n$ and taking expectation, we get
\begin{align*}
{\mathbb E}\left[{\cal A}_0 g\big(B_{k_1,k_2}^n\big)\right] &= \alpha \check{p}~ {\mathbb E}\big[g\big(B_{k_1,k_2}^n+1\big)\big]-\check{p}~{\mathbb E}\big[B_{k_1,k_2}^ng\big(B_{k_1,k_2}^n+1\big)\big]-\check{q}~{\mathbb E}\big[B_{k_1,k_2}^ng\big(B_{k_1,k_2}^n\big)\big]\\
&=\alpha \check{p}~ {\mathbb E}\big[g\big(B_{k_1,k_2}^n+1\big)\big]-{\mathbb E}\big[B_{k_1,k_2}^ng\big(B_{k_1,k_2}^n\big)\big]-\check{p}~{\mathbb E}\big[B_{k_1,k_2}^n \Delta g\big(B_{k_1,k_2}^n\big)\big]\\
&= \sum_{l=1}^{n-k_1-k_2+1}{\mathbb E}(I_l) {\mathbb E}\big[g\big(B_{k_1,k_2}^n+1\big)\big]-\sum_{l=1}^{n-k_1-k_2+1}{\mathbb E}\big[I_l g\big(W_l+1\big)\big]-\check{p}~{\mathbb E}\big[B_{k_1,k_2}^n \Delta g\big(B_{k_1,k_2}^n\big)\big]\\
&~~~~~~~~~~~~~~~~~~~~~~~~~~~~~~~~~~~~~~~~~~~~~~~~~~~~~~~~~~~~~~~~~~~~~~~~~~~~~~~~~~~~~~~~~~~~~~~~~~~~~~~~~\{{\rm from~\eqref{two:mean2}}\}\\
&= \sum_{l=1}^{n-k_1-k_2+1}{\mathbb E}(I_l) {\mathbb E}\big[g\big(B_{k_1,k_2}^n+1\big)-g(X_l+1)\big]-\sum_{l=1}^{n-k_1-k_2+1}{\mathbb E}\big[I_l \big(g(W_l+1)-g(X_l+1)\big)\big]\\
&~~~-\check{p}~{\mathbb E}\big[B_{k_1,k_2}^n \Delta g\big(B_{k_1,k_2}^n\big)\big]\\
&=\sum_{l=1}^{n-k_1-k_2+1}{\mathbb E}(I_l) \sum_{|u-l|\le k_1+k_2-1} {\mathbb E}\big[I_u \Delta g(X_l+D_{l,u-1}+1)\big]-\check{p}~{\mathbb E}\big[B_{k_1,k_2}^n \Delta g\big(B_{k_1,k_2}^n\big)\big]\\
&~~~-\sum_{l=1}^{n-k_1-k_2+1}\sum_{\substack{|u-l|\le k_1+k_2-1\\ u\neq l}}{\mathbb E}\big[I_l I_u \Delta g(X_l+\tilde{D}_{l,u-1}+1))\big].
\end{align*}
Note that
\begin{align*}
\sum_{l=1}^{n-k_1-k_2+1}\sum_{\substack{|u-l|\le k_1+k_2-1\\ u\neq l}}&{\mathbb E}\big[I_l I_u \Delta g(X_l+\tilde{D}_{l,u-1}+1))\big]\\
&=\sum_{l=1}^{n-k_1-k_2+1}\sum_{\substack{|u-l|\le k_1+k_2-1\\ u\neq l}}{\mathbb E}\big[I_l I_u\big] {\mathbb E}\big[ \Delta g(X_l+\tilde{D}_{l,u-1}+1)\big|I_l=I_u=1\big]=0.
\end{align*}
Therefore,
\begin{align*}
{\mathbb E}\left[{\cal A}_0 g\big(B_{k_1,k_2}^n\big)\right]&=\sum_{l=1}^{n-k_1-k_2+1}{\mathbb E}(I_l) \sum_{|u-l|\le k_1+k_2-1} {\mathbb E}\big[I_u \Delta g(X_l+D_{l,u-1}+1)\big]-\check{p}~{\mathbb E}\big[B_{k_1,k_2}^n \Delta g\big(B_{k_1,k_2}^n\big)\big].
\end{align*}
Hence, for $g \in {\cal G}_Z \cap {\cal G}_{B_{k_1,k_2}^n}$, we get
\begin{align*}
\big|{\mathbb E}\left[{\cal A}_0 g\big(B_{k_1,k_2}^n\big)\right]\big|\le \|\Delta g\|\sum_{l=1}^{n-k_1-k_2+1}a(p_l)\left\{\sum_{|u-l|\le k_1+k_2-1} a(p_u)+\check{p}\right\}.
\end{align*}
Using \eqref{two:bound}, the proof follows.\qed\\
Next, using an appropriate coupling, we derive two-parameter approximation result. For the ease of derivation, we define an auxiliary random variable $M$ as follows:\\
Let $\{\xi_l,1\le l \le n+k_1+k_2-1\}$ be independent Bernoulli trials with probability ${\mathbb P}(\xi_l=1)=p_l$. Define
\begin{equation}
M=\sum_{l=1}^{n}(1-\xi_{l})\dotsb (1-\xi_{l+k_1-1}) \xi_{l+k_1}\dotsb \xi_{l+k_1+k_2-1},\label{two:M}
\end{equation}
where $\xi_{l+mn}$ is treated as $\xi_l$, $1\le l \le n$ and $m=\pm1,\pm2,\dotsc$. From \eqref{two:M}, the mean and variance of $M$ are
\begin{align}
{\mathbb E}(M)=\sum_{l=1}^{n}a(p_l)\quad{\rm and}\quad Var(M)=\sum_{l=1}^{n}[a(p_l)-(a(p_l))^2]-2\sum_{\substack{l<r\\r-l\le k_1+k_2-1}}a(p_l)a(p_r).\label{two:MMV}
\end{align}
Next, it follows from \eqref{two:MMV} that ${\mathbb E}(M)>Var(M)$. Therefore, a pseudo-binomial approximation is suitable for $M$. In order to study $Z$-approximation to the distribution of $M$, let
\vspace{-0.07cm}
\begin{equation}
\alpha \check{p}={\mathbb E}(M)\quad\quad{\rm and}\quad\quad\alpha \check{p}\check{q}=Var(M).\label{two:MMAV}
\end{equation}
\vspace{-0.08cm}
Next, define
\vspace{-0.5cm}
\begin{align*}
M_l&=M-I_l, \quad\quad\quad\quad\quad\quad\quad\quad\quad\quad\quad\quad N_l=M-\sum_{|u-l|\le k_1+k_2-1}I_u\\
Q_l&=M-\sum_{|u-l|\le 2(k_1+k_2-1)}I_u\quad\quad{\rm and}\quad\quad S_l=M-\sum_{|u-l|\le 3(k_1+k_2-1)}I_u.
\end{align*}
Then, $N_l$ and $I_l$ are independent, $Q_l$ and $\{I_u:|u-l|\le k_1+k_2-1\}$ are independent, and $S_l$ and $\{I_u:|u-l|\le 2(k_1+k_2-1)\}$ are independent. Also, let $V_m$ be the $m$-th largest number of $(1-p_{l})\dotsb(1-p_{l+k_1-1})p_{l+k_1}\dotsb p_{l+k_1+k_2-1}(p_{l-2}^2(1-p_{l-1})p_l+p_{l+1}p_{l-1}^2)$, $2 \le l \le n$ and
$$\Psi:=2 \wedge \frac{4.6}{\sqrt{\sum_{m=4(k_1+k_2)-1}^{n}V_m}}.$$
We can now generalize Lemma $2.1$ of Wang and Xia \cite{WX} to our setting as follows.
\begin{lemma}\label{two:cople}
Let $\xi_0=1$ and $\xi_1,\xi_2,\dotsc,\xi_{n+k_1+k_2-1}$ be independent Bernoulli random variables with ${\mathbb P}(\xi_l=1)=p_l$. Then, for $n\ge k_1+k_2$
\begin{align*}
d_{TV}(B_{k_1,k_2}^{n+k_1+k_2-1}&,B_{k_1,k_2}^{n+k_1+k_2-1}+1)\\
&\le 1 \wedge \frac{2.3}{\sqrt{\sum_{l=2}^{n}(1-p_{l-1})\dotsb(1-p_{l+k_1-1})p_{l+k_1}\dotsb p_{l+k_1+k_2-1}(p_{l-2}^2(1-p_{l-1})p_l+p_{l+1}p_{l-1}^2)}}\\
&:=C(p_1,p_2,\dotsc,p_{n+k_1+k_2-1}).
\end{align*}
\end{lemma}
\noindent
{\bf Proof}. Let $B_{k_1,k_2}^m=\sum_{l=1}^{m}(1-\xi_l)\dotsb(1-\xi_{l+k_1-1})\xi_{l+k_1}\dotsb\xi_{l+k_1+k_2-1}$ and $\zeta_l$ be an independent copy of $\xi_l$ for $1\le l \le n+k_1+k_2-1$. Define $\hat{\xi}_0=1$ and, for $1\le l \le n+k_1+k_2-1$,
$$\hat{\xi}_l=\left\{
          \begin{array}{ll}
          \zeta_l & \text{if }\xi_{l-1}=\hat{\xi}_{l-1}=0=1-\xi_{l-2}=1-\hat{\xi}_{l-2},\\
          \xi_l & \text{otherwise},
          \end{array}\right.
$$
Let $\hat{B}_{k_1,k_2}^m=1+\sum_{l=1}^{m}(1-\hat{\xi}_l)\dotsb(1-\hat{\xi}_{l+k_1-1})\hat{\xi}_{l+k_1}\dotsb\hat{\xi}_{l+k_1+k_2-1}$ and
\begin{align*}
E_m&~=B_{k_1,k_2}^m - \hat{B}_{k_1,k_2}^m:=\sum_{l=1}^{m}\delta_l-1\\
&~=\sum_{l=1}^{m}\big[(1-\xi_l)\dotsb(1-\xi_{l+k_1-1})\xi_{l+k_1}\dotsb\xi_{l+k_1+k_2-1}-(1-\hat{\xi}_l)\dotsb(1-\hat{\xi}_{l+k_1-1})\hat{\xi}_{l+k_1}\dotsb\hat{\xi}_{l+k_1+k_2-1}\big]-1
\end{align*}
with $E_0=-1$. Now, observe that $\delta_l$ can take the values $0$ and $\pm1$, and $\{E_m, 1\le m \le n\}$ is a symmetric random walk. Also,
\begin{align*}
\{\delta_l\hspace{-0.07cm}=\hspace{-0.07cm}1\}&\hspace{-0.07cm}=\hspace{-0.07cm}\{\xi_{l-1}\hspace{-0.07cm}=\hspace{-0.07cm}\hat{\xi}_{l-1}\hspace{-0.07cm}=\hspace{-0.07cm}0\hspace{-0.07cm}=\hspace{-0.07cm}1\hspace{-0.07cm}-\hspace{-0.07cm}\xi_{l-2}\hspace{-0.07cm}=\hspace{-0.07cm}1\hspace{-0.07cm}-\hspace{-0.07cm}\hat{\xi}_{l-2},\hat{\xi}_l-{\xi}_l\hspace{-0.07cm}=\hspace{-0.07cm}1,\xi_{l+1}=\dotsb=\xi_{l+k_1-1}\hspace{-0.07cm}=\hspace{-0.07cm}0,\xi_{l+k_1}\hspace{-0.07cm}=\hspace{-0.07cm}\dotsb\hspace{-0.07cm}=\hspace{-0.07cm}\xi_{l+k_1+k_2-1}\hspace{-0.07cm}=\hspace{-0.07cm}1\}\\
&~~~\cup \{\xi_{l}\hspace{-0.07cm}=\hspace{-0.07cm}\hat{\xi}_{l}\hspace{-0.07cm}=\hspace{-0.07cm}0\hspace{-0.07cm}=\hspace{-0.07cm}1\hspace{-0.07cm}-\hspace{-0.07cm}\xi_{l-1}\hspace{-0.07cm}=\hspace{-0.07cm}1\hspace{-0.07cm}-\hspace{-0.07cm}\hat{\xi}_{l-1},\hat{\xi}_{l+1}\hspace{-0.07cm}-\hspace{-0.07cm}{\xi}_{l+1}\hspace{-0.07cm}=\hspace{-0.07cm}1,\xi_{l+2}\hspace{-0.07cm}=\hspace{-0.07cm}\dotsb\hspace{-0.07cm}=\hspace{-0.07cm}\xi_{l+k_1-1}\hspace{-0.07cm}=\hspace{-0.07cm}0,\xi_{l+k_1}\hspace{-0.07cm}=\hspace{-0.07cm}\dotsc\hspace{-0.07cm}=\hspace{-0.07cm}\xi_{l+k_1+k_2-1}\hspace{-0.07cm}=\hspace{-0.07cm}1\},\\
\{\delta_l\hspace{-0.07cm}=\hspace{-0.07cm}-1\}&\hspace{-0.07cm}=\hspace{-0.07cm}\{\xi_{l-1}\hspace{-0.07cm}=\hspace{-0.07cm}\hat{\xi}_{l-1}\hspace{-0.07cm}=\hspace{-0.07cm}0\hspace{-0.07cm}=\hspace{-0.07cm}1\hspace{-0.07cm}-\hspace{-0.07cm}\xi_{l-2}\hspace{-0.07cm}=\hspace{-0.07cm}1\hspace{-0.07cm}-\hspace{-0.07cm}\hat{\xi}_{l-2},\hat{\xi}_l-\xi_l\hspace{-0.07cm}=\hspace{-0.07cm}-\hspace{-0.07cm}1,\xi_{l+1}\hspace{-0.07cm}=\hspace{-0.07cm}\dotsb\hspace{-0.07cm}=\hspace{-0.07cm}\xi_{l+k_1-1}\hspace{-0.07cm}=\hspace{-0.07cm}0,\xi_{l+k_1}\hspace{-0.07cm}=\hspace{-0.07cm}\dotsc\hspace{-0.07cm}=\hspace{-0.07cm}\xi_{l+k_1+k_2-1}\hspace{-0.07cm}=\hspace{-0.07cm}1\}\\
&~~~\cup \{\xi_{l}\hspace{-0.07cm}=\hspace{-0.07cm}\hat{\xi}_{l}\hspace{-0.07cm}=\hspace{-0.07cm}0\hspace{-0.07cm}=\hspace{-0.07cm}1\hspace{-0.07cm}-\hspace{-0.07cm}\xi_{l-1}\hspace{-0.07cm}=\hspace{-0.07cm}1\hspace{-0.07cm}-\hspace{-0.07cm}\hat{\xi}_{l-1},\hat{\xi}_{l+1}\hspace{-0.07cm}-\hspace{-0.07cm}\xi_{l+1}\hspace{-0.07cm}=\hspace{-0.07cm}-1,\xi_{l+2}\hspace{-0.07cm}=\hspace{-0.07cm}\dotsb\hspace{-0.07cm}=\hspace{-0.07cm}\xi_{l+k_1-1}\hspace{-0.07cm}=\hspace{-0.07cm}0,\xi_{l+k_1}\hspace{-0.07cm}=\hspace{-0.07cm}\dotsc\hspace{-0.07cm}=\hspace{-0.07cm}\xi_{l+k_1+k_2-1}\hspace{-0.07cm}=\hspace{-0.07cm}1\}
\end{align*}
and hence
\begin{align*}
\{\delta_l\neq0\}&\hspace{-0.07cm}=\hspace{-0.07cm}\{\xi_{l-1}\hspace{-0.07cm}=\hspace{-0.07cm}\hat{\xi}_{l-1}\hspace{-0.07cm}=\hspace{-0.07cm}0\hspace{-0.07cm}=\hspace{-0.07cm}1\hspace{-0.07cm}-\hspace{-0.07cm}\xi_{l-2}\hspace{-0.07cm}=\hspace{-0.07cm}1\hspace{-0.07cm}-\hspace{-0.07cm}\hat{\xi}_{l-2},|\hat{\xi}_l\hspace{-0.07cm}-\hspace{-0.07cm}{\xi}_l|\hspace{-0.07cm}=\hspace{-0.07cm}1,\xi_{l+1}=\dotsb=\xi_{l+k_1-1}\hspace{-0.07cm}=\hspace{-0.07cm}0,\xi_{l+k_1}\hspace{-0.07cm}=\hspace{-0.07cm}\dotsc\hspace{-0.07cm}=\hspace{-0.07cm}\xi_{l+k_1+k_2-1}\hspace{-0.07cm}=\hspace{-0.07cm}1\}\\
&~~~\cup \{\xi_{l}\hspace{-0.07cm}=\hspace{-0.07cm}\hat{\xi}_{l}\hspace{-0.07cm}=\hspace{-0.07cm}0\hspace{-0.07cm}=\hspace{-0.07cm}1\hspace{-0.07cm}-\hspace{-0.07cm}\xi_{l-1}\hspace{-0.07cm}=\hspace{-0.07cm}1-\hat{\xi}_{l-1},|\hat{\xi}_{l+1}\hspace{-0.07cm}-\hspace{-0.07cm}{\xi}_{l+1}|\hspace{-0.07cm}=\hspace{-0.07cm}1,\xi_{l+2}\hspace{-0.07cm}=\hspace{-0.07cm}\dotsb\hspace{-0.07cm}=\hspace{-0.07cm}\xi_{l+k_1-1}\hspace{-0.07cm}=\hspace{-0.07cm}0,\xi_{l+k_1}\hspace{-0.07cm}=\hspace{-0.07cm}\dotsc\hspace{-0.07cm}=\hspace{-0.07cm}\xi_{l+k_1+k_2-1}\hspace{-0.07cm}=\hspace{-0.07cm}1\}.
\end{align*}
Define $R_l={\bf 1}_{\{\delta_l\neq0\}}$ and $R=\sum_{l=2}^{n}R_l$. Therefore,
\begin{align}
{\mathbb E}(R_l)&={\mathbb P}(\delta_l=1)+{\mathbb P}(\delta_l=-1)=2{\mathbb P}(\delta_l=1)\nonumber\\
&=2\big[{\mathbb P}(\xi_{l-1}\hspace{-0.07cm}=\hspace{-0.07cm}\hat{\xi}_{l-1}\hspace{-0.07cm}=\hspace{-0.07cm}0\hspace{-0.07cm}=\hspace{-0.07cm}1\hspace{-0.07cm}-\hspace{-0.07cm}\xi_{l-2}\hspace{-0.07cm}=\hspace{-0.07cm}1\hspace{-0.07cm}-\hspace{-0.07cm}\hat{\xi}_{l-2},\zeta_l\hspace{-0.07cm}=\hspace{-0.07cm}1,\xi_l=\dotsb=\xi_{l+k_1-1}\hspace{-0.07cm}=\hspace{-0.07cm}0,\xi_{l+k_1}\hspace{-0.07cm}=\hspace{-0.07cm}\dotsc\hspace{-0.07cm}=\hspace{-0.07cm}\xi_{l+k_1+k_2-1}\hspace{-0.07cm}=\hspace{-0.07cm}1)\nonumber\\
&~~~+{\mathbb P}(\xi_{l}\hspace{-0.07cm}=\hspace{-0.07cm}\hat{\xi}_{l}\hspace{-0.07cm}=\hspace{-0.07cm}0\hspace{-0.07cm}=\hspace{-0.07cm}1\hspace{-0.07cm}-\hspace{-0.07cm}\xi_{l-1}\hspace{-0.07cm}=\hspace{-0.07cm}1\hspace{-0.07cm}-\hspace{-0.07cm}\hat{\xi}_{l-1},\zeta_{l+1}\hspace{-0.07cm}=\hspace{-0.07cm}1,\xi_{l+1}=\dotsb=\xi_{l+k_1-1}\hspace{-0.07cm}=\hspace{-0.07cm}0,\xi_{l+k_1}\hspace{-0.07cm}=\hspace{-0.07cm}\dotsc\hspace{-0.07cm}=\hspace{-0.07cm}\xi_{l+k_1+k_2-1}\hspace{-0.07cm}=\hspace{-0.07cm}1)\big]\nonumber\\
&=2\big[(1-p_l)p_l {\mathbb P}(\xi_{l-1}=\hat{\xi}_{l-1}=0=1-\xi_{l-2}=1-\hat{\xi}_{l-2})+p_{l+1}{\mathbb P}(\xi_{l}=\hat{\xi}_{l}=0=1-\xi_{l-1}=1-\hat{\xi}_{l-1})\big]\times\nonumber\\
&~~~~~~(1-p_{l+1})\dotsb(1-p_{l+k_1-1})p_{l+k_1}\dotsb p_{l+k_1+k_2-1}.\label{two:expe}
\end{align}
Now, consider
\begin{align}
&{\mathbb P}(\xi_{l-1}=\hat{\xi}_{l-1}=0=1-\xi_{l-2}=1-\hat{\xi}_{l-2})\nonumber\\
&={\mathbb P}(\xi_{l-1}=\hat{\xi}_{l-1}=0|\xi_{l-2}=\hat{\xi}_{l-2}=1){\mathbb P}(\xi_{l-2}=\hat{\xi}_{l-2}=1)\nonumber\\
&={\mathbb P}(\xi_{l-1}=0){\mathbb P}(\xi_{l-2}=\hat{\xi}_{l-2}=1)\nonumber\\
&=(1\hspace{-0.05cm}-\hspace{-0.05cm}p_{l-1})\big[{\mathbb P}(\xi_{l-2}=\hat{\xi}_{l-2}=1|\xi_{l-3}=\hat{\xi}_{l-3}=0=1\hspace{-0.05cm}-\hspace{-0.05cm}\xi_{l-4}=1\hspace{-0.05cm}-\hspace{-0.05cm}\hat{\xi}_{l-4}){\mathbb P}(\xi_{l-3}=\hat{\xi}_{l-3}=0=1\hspace{-0.05cm}-\hspace{-0.05cm}\xi_{l-4}=1\hspace{-0.05cm}-\hspace{-0.05cm}\hat{\xi}_{l-4})\nonumber\\
&~~~+{\mathbb P}(\xi_{l-2}=\hat{\xi}_{l-2}=1|\{\xi_{l-3}=\hat{\xi}_{l-3}=0=1-\xi_{l-4}=1-\hat{\xi}_{l-4}\}^{c}){\mathbb P}(\{\xi_{l-3}=\hat{\xi}_{l-3}=0=1-\xi_{l-4}=1-\hat{\xi}_{l-4}\}^c)\big]\nonumber\\
&=(1-p_{l-1})\big[{\mathbb P}(\xi_{l-2}=1){\mathbb P}(\zeta_{l-2}=1){\mathbb P}(\xi_{l-3}=\hat{\xi}_{l-3}=0=1-\xi_{l-4}=1-\hat{\xi}_{l-4})\nonumber\\
&~~~+{\mathbb P}(\xi_{l-2}=1){\mathbb P}(\{\xi_{l-3}=\hat{\xi}_{l-3}=0=1-\xi_{l-4}=1-\hat{\xi}_{l-4}\}^c)\big]\ge (1-p_{l-1})p_{l-2}^2,\label{two:gg1}
\end{align}
where $A^c$ denotes the complement of set $A$. Similarly,
\begin{equation}
{\mathbb P}(\xi_{l}=\hat{\xi}_{l}=0=1-\xi_{l-1}=1-\hat{\xi}_{l-1}) \ge  (1-p_{l})p_{l-1}^2.\label{two:gg2}
\end{equation}
Substituting \eqref{two:gg1} and \eqref{two:gg2} in \eqref{two:expe}, we get
$${\mathbb E}(R_l) \ge 2(1-p_{l})\dotsb(1-p_{l+k_1-1})p_{l+k_1}\dotsb p_{l+k_1+k_2-1}(p_{l-2}^2(1-p_{l-1})p_l+p_{l+1}p_{l-1}^2).$$
Hence,
\begin{equation}
{\mathbb E}(R) \ge 2\sum_{l=2}^{n}(1-p_{l})\dotsb(1-p_{l+k_1-1})p_{l+k_1}\dotsb p_{l+k_1+k_2-1}(p_{l-2}^2(1-p_{l-1})p_l+p_{l+1}p_{l-1}^2).\label{two:expR}
\end{equation}
Next, we calculate the upper bound for variance of $R$. Observe that
$$Var(R)=\sum_{l=2}^{n}({\mathbb E}(R_l)-({\mathbb E}(R_l))^2)+2\sum_{l<r}Cov(R_l,R_r)\le {\mathbb E}(R)+2\sum_{l=2}^{n}\sum_{r=l+k_1+k_2}^{n}[{\mathbb E}(R_lR_r)-{\mathbb E}(R_l){\mathbb E}(R_r)].$$
Since $R_lR_r=0$ for $r-l\le k_1+k_2-1$. Now, for $r-l=k_1+k_2$ and $k_1+k_2+1$, it can be verified that
\begin{align*}
|{\mathbb E}(R_lR_r)-{\mathbb E}(R_l){\mathbb E}(R_r)| &\le{\mathbb E}(R_l)|{\mathbb E}(R_r|\xi_{l+k_1+k_2-1}=\hat{\xi}_{l+k_1+k_2-1}=1)-{\mathbb E}(R_r|\xi_{l+k_1+k_2-1}=\hat{\xi}_{l+k_1+k_2-1}=0)|.
\end{align*}
For $r-l \ge k_1+k_2+2$, consider
\begin{align*}
{\mathbb E}(R_l R_r) = {\mathbb P}(R_l=1,R_r=1)= {\mathbb P}(R_r=1|R_l=1){\mathbb P}(R_l=1)={\mathbb E}(R_l) {\mathbb E}(R_r|\xi_{l+k_1+k_2-1}=\hat{\xi}_{l+k_1+k_2-1}=1).
\end{align*}
Note that ${\mathbb E}(R_r|\xi_{l+k_1+k_2-1}=\hat{\xi}_{l+k_1+k_2-1}=1)={\mathbb E}(R_r|\xi_{l+k_1+k_2-1}=0,\hat{\xi}_{l+k_1+k_2-1}=1)={\mathbb E}(R_r|\xi_{l+k_1+k_2-1}=1,\hat{\xi}_{l+k_1+k_2-1}=0)$. Therefore, for $r-l \ge k_1+k_2$
\begin{align*}
|{\mathbb E}(R_lR_r)-{\mathbb E}(R_l){\mathbb E}(R_r)| &= {\mathbb E}(R_l)|{\mathbb E}(R_r|\xi_{l+k_1+k_2-1}=\hat{\xi}_{l+k_1+k_2-1}=1)-{\mathbb E}R_r|\\
&\le {\mathbb E}R_l|{\mathbb E}(R_r|\xi_{l+k_1+k_2-1}=\hat{\xi}_{l+k_1+k_2-1}=1)-{\mathbb E}(R_r|\xi_{l+k_1+k_2-1}=\hat{\xi}_{l+k_1+k_2-1}=0)|.
\end{align*}
Let $W=\min\{m\ge l+k_1+k_2:\xi_m=\zeta_m\}$. Then $W$ is independent of $\xi_l$ and $\hat{\xi_l}$ with ${\mathbb P}(W\ge r+1)= 2\prod_{m=l+k_1+k_2}^{r}p_m(1-p_m)$ and
$$((\xi_i,\hat{\xi}_i)_{i\ge s}| W=s, \xi_{l+k_1+k_2-1}=\hat{\xi}_{l+k_1+k_2-1}=1)\stackrel{\cal L}{=}((\xi_i,\hat{\xi}_i)_{i\ge s}| W=s, \xi_{l+k_1+k_2-1}=\hat{\xi}_{l+k_1+k_2-1}=0)$$
for all $l+ k_1+k_2 \le s \le r$. Hence,
$${\mathbb E}(R_r|W=m,\xi_{l+k_1+k_2-1}=\hat{\xi}_{l+k_1+k_2-1}=1)={\mathbb E}(R_r|W=m,\xi_{l+k_1+k_2-1}=\hat{\xi}_{l+k_1+k_2-1}=0).$$
Since $2p_m(1-p_m)\le1/2$, Therefore,
\begin{align*}
|{\mathbb E}(R_lR_r)-{\mathbb E}(R_l){\mathbb E}(R_r)| &= {\mathbb E}(R_l){\mathbb P}(W\ge r+1)\le {\mathbb E}(R_l) 2^{-(r-l-k_1-k_2+1)}.
\end{align*}
Therefore,
\begin{align}
Var(R)\le {\mathbb E}(R)+2\sum_{l=2}^{n}\sum_{r=l+k_1+k_2}^{n}{\mathbb E}(R_l) 2^{-(r-l-k_1-k_2+1)}\le3 {\mathbb E}(R).\label{two:vvar}
\end{align}
The remaining part of the proof follows from Lemma $2.1$ of Wang and Xia \cite{WX}. Define $\tau_0=0$ and
$$\tau_l=\min\{m>\tau_{l-1}:\delta_m\neq0\}.$$
So, $R=r$ implies $\tau_r \le n < \tau_{r+1}$ and the distribution of the random variable $(\delta_{\tau_l}|R=r)$, $1\le l \le r$ is uniformly distributed on $\{-1,1\}$. Let $W_l^r=(E_{\tau_l}|R=r)$, $1\le l \le r$. Then, $W_l^r$ is a simple random walk with initial state $-1$, i.e.
$${\mathbb P}(W_l^r-W_{l-1}^r=1)={\mathbb P}(W_l^r-W_{l-1}^r=-1)=1/2.$$
Let $K=\min\{l\ge1:E_{\tau_l}=0\}$ and
$$\check{B}_{k_1,k_2}^m=\left\{
          \begin{array}{ll}
          \hat{B}_{k_1,k_2}^m & \text{if }m < \tau_K,\\
          B_{k_1,k_2}^m & \text{if } m\ge \tau_K.
          \end{array}\right.
$$
Therefore, by coupling, $\check{B}_{k_1,k_2}^{n+k_1,k_2-1}\stackrel{\cal L}{=}B_{k_1,k_2}^{n+k_1+k_2-1}+1 $ and
$$d_{TV}(B_{k_1,k_2}^{n\hspace{-0.01cm}+\hspace{-0.01cm}k_1\hspace{-0.01cm}+\hspace{-0.01cm}k_2\hspace{-0.01cm}-\hspace{-0.01cm}1},B_{k_1,k_2}^{n\hspace{-0.01cm}+\hspace{-0.01cm}k_1\hspace{-0.01cm}+\hspace{-0.01cm}k_2\hspace{-0.01cm}-\hspace{-0.01cm}1}\hspace{-0.07cm}+\hspace{-0.07cm}1) \le d_{TV}(B_{k_1,k_2}^{n\hspace{-0.01cm}+\hspace{-0.01cm}k_1\hspace{-0.01cm}+\hspace{-0.01cm}k_2\hspace{-0.01cm}-\hspace{-0.01cm}1},\check{B}_{k_1,k_2}^{n\hspace{-0.01cm}+\hspace{-0.01cm}k_1\hspace{-0.01cm}+\hspace{-0.01cm}k_2\hspace{-0.01cm}-\hspace{-0.01cm}1}\hspace{-0.07cm})\le {\mathbb P}\big(B_{k_1,k_2}^{n\hspace{-0.01cm}+\hspace{-0.01cm}k_1\hspace{-0.01cm}+\hspace{-0.01cm}k_2\hspace{-0.01cm}-\hspace{-0.01cm}1}\neq\check{B}_{k_1,k_2}^{n\hspace{-0.01cm}+\hspace{-0.01cm}k_1\hspace{-0.01cm}+\hspace{-0.01cm}k_2\hspace{-0.01cm}-\hspace{-0.01cm}1}\hspace{-0.07cm}\big).$$
By the reflection principle, it can be easily seen that
\begin{align*}
{\mathbb P}\big(B_{k_1,k_2}^{n+k_1+k_2-1}\neq\check{B}_{k_1,k_2}^{n+k_1+k_2-1}|R=r\big)&={\mathbb P}\left(\max_{1\le l\le r}E_{\tau_l}\le-1|R=r\right)={\mathbb P}\left(\max_{1\le l\le r}W_l^r\le -1\right)~~~~~~~~~~~~~~~~~~~~
\end{align*}
\begin{align*}
&~~~~~~~~~~~~~~~~~~~~~~~~~~~~~~~~~~~={\mathbb P}\big(W_r^r\in\{-2,-1\}\big)=\max_{m}{\mathbb P}\big(W_r^r=m\big)\le\sqrt{\frac{2}{\pi r}}\le \frac{0.8}{\sqrt{r}}.
\end{align*}
Therefore, for all $r \ge 1$ and $0 < \kappa < 1$,
\begin{align}
{\mathbb P}\big(B_{k_1,k_2}^{n+k_1+k_2-1}\neq\check{B}_{k_1,k_2}^{n+k_1+k_2-1}\big)&=\sum_{r=1}^{\infty}{\mathbb P}\big(B_{k_1,k_2}^{n+k_1+k_2-1}\neq\check{B}_{k_1,k_2}^{n+k_1+k_2-1}|R=r\big){\mathbb P}(R=r)\nonumber\\
&\le \sum_{1\le r\le \kappa {\mathbb E}(R)}{\mathbb P}(R=r)+\sum_{r > \kappa {\mathbb E}(R)} \frac{0.8}{\sqrt{\kappa {\mathbb E}(R)}}{\mathbb P}(R=r)\nonumber\\
&\le \frac{Var(R)}{(1-\kappa)^2({\mathbb E}(R))^2}+\frac{0.8}{\sqrt{\kappa {\mathbb E}(R)}}. \label{two:eer}
\end{align}
where the last inequality is obtained by Chebyshev's inequality.\\
Now, if $\sum_{l=2}^{n}(1-p_{l})\dotsb(1-p_{l+k_1-1})p_{l+k_1}\dotsb p_{l+k_1+k_2-1}(p_{l-2}^2(1-p_{l-1})p_l+p_{l+1}p_{l-1}^2)\le 5.29$ then proof is follows. So, let us consider $\sum_{l=2}^{n}(1-p_{l})\dotsb(1-p_{l+k_1-1})p_{l+k_1}\dotsb p_{l+k_1+k_2-1}(p_{l-2}^2(1-p_{l-1})p_l+p_{l+1}p_{l-1}^2)>5.29$. Then ${\mathbb E}(R)>10.58$. Hence, using \eqref{two:expR},\eqref{two:vvar} and \eqref{two:eer}, we get
\begin{align*}
{\mathbb P}\big(B_{k_1,k_2}^{n\hspace{-0.01cm}+\hspace{-0.01cm}k_1\hspace{-0.01cm}+\hspace{-0.01cm}k_2\hspace{-0.01cm}-\hspace{-0.01cm}1}\neq\check{B}_{k_1,k_2}^{n\hspace{-0.01cm}+\hspace{-0.01cm}k_1\hspace{-0.01cm}+\hspace{-0.01cm}k_2\hspace{-0.01cm}-\hspace{-0.01cm}1}\hspace{-0.07cm}\big)&\le \frac{3}{(1-\kappa)^2{\mathbb E}(R)}+\frac{0.8}{\sqrt{\kappa {\mathbb E}(R)}}\\
&\le  \hspace{-0.08cm}\frac{2.3}{\sqrt{\sum_{l=2}^{n}(1\hspace{-0.07cm}-\hspace{-0.07cm}p_{l\hspace{-0.01cm}-\hspace{-0.01cm}1})\dotsb(1\hspace{-0.07cm}-\hspace{-0.07cm}p_{l\hspace{-0.01cm}+\hspace{-0.01cm}k_1\hspace{-0.01cm}-\hspace{-0.01cm}1})p_{l\hspace{-0.01cm}+\hspace{-0.01cm}k_1}\dotsb p_{l\hspace{-0.01cm}+\hspace{-0.01cm}k_1\hspace{-0.01cm}+\hspace{-0.01cm}k_2\hspace{-0.01cm}-\hspace{-0.01cm}1}(p_{l\hspace{-0.01cm}-\hspace{-0.01cm}2}^2(1\hspace{-0.07cm}-\hspace{-0.07cm}p_{l\hspace{-0.01cm}-\hspace{-0.01cm}1})p_l\hspace{-0.07cm}+\hspace{-0.07cm}p_{l\hspace{-0.01cm}+\hspace{-0.01cm}1}p_{l\hspace{-0.01cm}-\hspace{-0.01cm}1}^2)}},
\end{align*}
where the last inequality is obtained by assuming $\kappa=0.2197$. This proves the lemma.\qed

\begin{theorem}\label{two:Ith2}
Let $n\ge 4(k_1+k_2)$ and $M$ be as defined in \eqref{two:M} with \eqref{two:MMAV}. Then
\begin{align}
d_{TV}(M,Z)&\le\frac{\Psi}{\floor{\alpha}\check{p}\check{q}}\sum_{l=1}^{n}a(p_l)\left\{\sum_{|u-l|\le k_1+k_2-1}\sum_{v=l-2(k_1+k_2)+2}^{u-k_1-k_2}a(p_u)a(p_v)+\sum_{u=l-k_1-k_2+1}^{l}\sum_{v=l+k_1+k_2}^{l+2(k_1+k_2-1)}a(p_u)a(p_v)\right.\nonumber\\
&~~~\left.+\sum_{u=l+1}^{l+k_1+k_2-1}\sum_{v=u+k_1+k_2}^{u+2(k_1+k_2-1)}\hspace{-0.33cm}a(p_u)a(p_v)+\sum_{|u-l|\le k_1+k_2-1}\sum_{|v-l|\le2(k_1+k_2-1)}\hspace{-0.33cm}a(p_u)a(p_v)+\check{p}\sum_{|u-l|\le k_1+k_2-1}\hspace{-0.33cm}a(p_u)\right\}\label{two:eee1}\\
&\le \frac{\Psi}{\floor{\alpha}\check{p}\check{q}}\sum_{l=1}^{n}a(p_l)\left\{2\sum_{|u-l|\le k_1+k_2-1}\sum_{|v-l|\le2(k_1+k_2-1)}a(p_u)a(p_v)+\check{p}\sum_{|u-l|\le k_1+k_2-1}a(p_u)\right\}\label{two:eee2}\\
&\le \frac{(2k_1+2k_2-1)\Psi}{\floor{\alpha}\check{p}\check{q}}\sum_{i=1}^{n}a(p_l)a(m_l)\left\{2(4k_1+4 k_2-3)a(m_l)+\check{p}\right\}\label{two:eee3},
\end{align}
where $m_l=\max\{p_s:|l-s|\le2(k_1+k_2-1)\}$.
\end{theorem}
\noindent
{\bf Proof}. Using the steps similar to the proof of Theorem \ref{two:Ith1} for the random variable  $M$, we get
$${\mathbb E}[{\cal A}_0 g(M)]=\sum_{l=1}^{n}{\mathbb E}(I_l)\sum_{|u-l|\le k_1+k_2-1}{\mathbb E}(I_u\Delta g(N_l+D_{l,u-1}+1))-\check{p}\sum_{l=1}^{n}{\mathbb E}(I_l \Delta g(M_l+1)).$$
From \eqref{two:MMAV}, it can be seen that
$$\check{p}\sum_{l=1}^{n}{\mathbb E}(I_l)-\sum_{l=1}^{n}{\mathbb E}(I_l)\sum_{|u-l|\le k_1+k_2-1}{\mathbb E}(I_u)=0.$$
Now, note that
\begin{align}
{\mathbb E}[{\cal A}_0 g(M)]&=\sum_{l=1}^{n}{\mathbb E}(I_l)\sum_{|u-l|\le k_1+k_2-1}{\mathbb E}(I_u\Delta g(N_l+D_{l,u-1}+1))-\check{p}\sum_{l=1}^{n}{\mathbb E}(I_l \Delta g(M_l+1))\nonumber\\
&~~~+\left[\check{p}\sum_{l=1}^{n}{\mathbb E}(I_l)-\sum_{l=1}^{n}{\mathbb E}(I_l)\sum_{|u-l|\le k_1+k_2-1}{\mathbb E}(I_u)\right]{\mathbb E}(\Delta g(M+1))\nonumber\\
&=\sum_{l=1}^{n}{\mathbb E}(I_l)\sum_{|u-l|\le k_1+k_2-1}[{\mathbb E}(I_u\Delta g(N_l+D_{l,u-1}+1))-{\mathbb E}(I_u) {\mathbb E}(\Delta g(M+1))]\nonumber\\
&~~~-\check{p}\sum_{l=1}^{n}[{\mathbb E}(I_l \Delta g(M_l+1))-{\mathbb E}(I_l){\mathbb E}(\Delta g(M+1))].\label{two:mainexp}
\end{align}
The following expression, from the first term of \eqref{two:mainexp}, can be rewritten as
\begin{align}
&{\mathbb E}(I_u\Delta g(N_l+D_{l,u-1}+1))-{\mathbb E}(I_u) {\mathbb E}(\Delta g(M+1))\nonumber\\
&={\mathbb E}[(I_u(\Delta g(N_l+D_{l,u-1}+1)-\Delta g(Q_l+1))]+{\mathbb E}(I_u){\mathbb E}[\Delta g(Q_l+1) -{\mathbb E}(\Delta g(M+1))].\label{two:hh}
\end{align}
Since $Q_l$ is independent of $I_u$. Further, to simplify calculations, observe that
\begin{align}
{\mathbb E}&[(I_u(\Delta g(N_l+D_{l,u-1}+1)- \Delta g(Q_l+1))]\nonumber\\
&=\sum_{v=l+k_1+k_2}^{l+2(k_1+k_2-1)}{\mathbb E}\left(I_u I_v \Delta^2 g\left(Q_l+\sum_{s=l-2(k_1+k_2-1)}^{u-1}I_s+\sum_{s=l+k_1+k_2}^{v-1}I_s+1\right)\right)\nonumber\\
&~~~+\sum_{v=l-2(k_1+k_2-1)}^{u-1}{\mathbb E}\left(I_u I_v \Delta^2 g\left(Q_l+\sum_{s=l-2(k_1+k_2-1)}^{v-1}I_s+1\right)\right)\nonumber\\
&=\sum_{v=l+k_1+k_2}^{l+2(k_1+k_2-1)}{\mathbb E}(I_u I_v) {\mathbb E}\left(\Delta^2 g\left(S_l+\sum_{s=l-3(k_1+k_2-1)}^{u-1}I_s+\sum_{s=l+k_1+k_2}^{v-1}I_s+\sum_{s=l+2k_1+2k_2-1}^{l+3(k_1+k_2-1)}I_s+1\right)\Bigr| I_u=I_v=1\right)\nonumber\\
&~~~+\sum_{v=l-2(k_1+k_2-1)}^{u-1}{\mathbb E}(I_u I_v){\mathbb E}\left( \Delta^2 g\left(S_l+\sum_{s=l-3(k_1+k_2-1)}^{v-1}I_s+\sum_{s=l+2k_1+2k_2-1}^{l+3(k_1+k_2-1)}I_s+1\right)\Bigr| I_u=I_v=1\right).\label{two:ww1}
\end{align}
Also, let $b(\eta_1,\eta_2,\dotsc,\eta_m)=\sum_{l=1}^{m-k_1-k_2+1}(1-\eta_l)\dotsb(1-\eta_{l+k_1-1})\eta_{l+k_1}\dotsb \eta_{l+k_1+k_2-1}$, for a sequence of independent Bernoulli random variables $\eta_1,\eta_2,\dotsc,\eta_m$. Define
$${\bf (01)}_{s,m}=\underbrace{00\dotsc00}_{s}\underbrace{11\dotsc11}_{m}$$
and if $s\le0$ and $m \le0$ respectively, then ${\bf (01)}_{s,m}=\underbrace{11\dotsc11}_{m}\quad{\rm and }\quad {\bf (01)}_{s,m}=\underbrace{00\dotsc00}_{s}$ respectively. Also,  ${\bf (01)}_{s,m}$ is an empty if $s$ and $m\le0$. Then,
\begin{align}
&\left( \left(S_l+\sum_{s=l-3(k_1+k_2-1)}^{v-1}I_s+\sum_{s=l+2k_1+2k_2-1}^{l+3(k_1+k_2-1)}I_s+1\right)\Bigr| I_u=I_v=1\right)\nonumber\\
&~~~~~~~\stackrel{\cal L}{=}b(\xi_{l+2k_1+2k_2-1},\dotsc,\xi_n,\xi_1,\dotsc,\xi_{v-1},{\bf (01)}_{k_1,k_2-1}),~{\rm and}~~~~~~~~\label{two:law}\\
&\left(\left(S_l+\sum_{s=l-3(k_1+k_2-1)}^{u-1}I_s+\sum_{s=l+k_1+k_2}^{v-1}I_s+\sum_{s=l+2k_1+2k_2-1}^{l+3(k_1+k_2-1)}I_s+1\right)\Bigr| I_u=I_v=1\right)\nonumber\\
&~~~~~~\stackrel{\cal L}{=}{b({\bf (01)}_{0,1\wedge (u-l)\wedge (v-l-k_1-k_2)},\xi_{(l\vee u)+k_1+k_2},\dotsc,\xi_{v-1},{\bf (01)}_{k_1,k_2},\xi_{v+k_1+k_2},\dotsc,\xi_n,\xi_1,\dotsc,\xi_{u-1},{\bf (01)}_{k_1,k_2-1})}.\nonumber
\end{align}
Now, we obtain error bounds for the terms involving in \eqref{two:hh}. Applying Lemma \ref{two:cople} and using \eqref{two:law}, we have
\begin{align}
&\left|{\mathbb E}\left( \Delta^2 g\left(S_l+\sum_{s=l-3(k_1+k_2-1)}^{v-1}I_s+\sum_{s=l+2k_1+2k_2-1}^{l+3(k_1+k_2-1)}I_s+1\right)\Bigr| I_u=I_v=1\right)\right|\nonumber\\
&~~~~~~~~~~~~~~~~~~~=|{\mathbb E}(\Delta^2 g(b(\xi_{l+2k_1+2k_2-1},\dotsc,\xi_n,\xi_1,\dotsc,\xi_{u-1},{\bf (01)}_{k_1,k_2-1})))|\nonumber\\
&~~~~~~~~~~~~~~~~~~~\le 2 \|\Delta g\|d_{TV}(b(\xi_{l+2k_1+2k_2-1},\dotsc,\xi_n,\xi_1,\dotsc,\xi_{u-1},{\bf (01)}_{k_1,k_2-1}),\nonumber\\
&~~~~~~~~~~~~~~~~~~~~~~~~~~~~~~~~~~~~~~b(\xi_{l+2k_1+2k_2-1},\dotsc,\xi_n,\xi_1,\dotsc,\xi_{u-1},{\bf (01)}_{k_1,k_2-1})+1)\nonumber\\
&~~~~~~~~~~~~~~~~~~~\le 2 \|\Delta g\|C(p_{l+2k_1+2k_2-1},\dotsc,p_n,p_1,\dotsc,p_{u-1},{\bf (01)}_{k_1,k_2-1})\nonumber\\
&~~~~~~~~~~~~~~~~~~~\le \Psi\|\Delta g\|.\label{two:ww2}
\end{align}
Similarly,
\begin{align}
\left| {\mathbb E}\left(\Delta^2 g\left(S_l+\sum_{s=l-3(k_1+k_2-1)}^{u-1}I_s+\sum_{s=l+k_1+k_2}^{v-1}I_s+\sum_{s=l+2k_1+2k_2-1}^{l+3(k_1+k_2-1)}I_s+1\right)\Bigr| I_u=I_v=1\right)\right|\le \Psi \|\Delta g\|.\label{two:ww3}
\end{align}
Substituting \eqref{two:ww2} and \eqref{two:ww3} in \eqref{two:ww1}, we have
\begin{align}
|{\mathbb E}[(I_u(\Delta g(N_l+D_{l,u-1}+1)-\Delta g(Q_l+1))]|\le \Psi\|\Delta g\|\left\{\sum_{v=l+k_1+k_2}^{l+2(k_1+k_2-1)}+\sum_{v=l-2(k_1+k_2-1)}^{u-1}\right\}{\mathbb E}(I_u I_v).\label{two:ww4}
\end{align}
Now, following similar steps for the second term of \eqref{two:mainexp} and \eqref{two:hh}, we get
\begin{align}
|{\mathbb E}[\Delta g(Q_l+1) -{\mathbb E}(\Delta g(M+1))]|&\le\Psi\|\Delta g\|\sum_{|v-l|\le 2(k_1+k_2-1)}{\mathbb E}(I_v).\label{two:ww5}
\end{align}
and
\begin{align}
|{\mathbb E}(I_l \Delta g(M_l+1))&-{\mathbb E}(I_l){\mathbb E}(\Delta g(M+1))|\nonumber\\
&=|{\mathbb E}(I_l( \Delta g(M_l+1)-\Delta g(N_l+1))+{\mathbb E}(I_l){\mathbb E}(\Delta(N_l+1)-\Delta g(M+1))|\nonumber\\
&=|{\mathbb E}(I_l (\Delta g(M_l+1)-\Delta g(N_l+1))|+{\mathbb E}(I_l)|{\mathbb E}(\Delta(N_l+1)-\Delta g(M+1))|\nonumber\\
&\le\Psi\|\Delta g\|\left\{\sum_{\substack{|u-l|\le k_1+k_2-1\\u\neq l}}{\mathbb E}(I_l I_u)+\sum_{|u-l|\le k_1+k_2-1}{\mathbb E}(I_l){\mathbb E}(I_u)\right\}.\label{two:ww6}
\end{align}
Substituting \eqref{two:ww4} and \eqref{two:ww5} in \eqref{two:hh}, and combining \eqref{two:mainexp} and \eqref{two:ww6}, we get
\begin{align*}
|{\mathbb E}[{\cal A}_0 g(M)]| &\le \Psi\|\Delta g\|\sum_{l=1}^{n}{\mathbb E}(I_l)\sum_{|u-l|\le k_1+k_2-1}\left\{\left(\sum_{v=l+k_1+k_2}^{l+2(k_1+k_2-1)}+\sum_{v=l-2(k_1+k_2-1)}^{u-1}\right){\mathbb E}(I_u I_v)+\sum_{|v-l|\le 2(k_1+k_2-1)}\hspace{-0.5cm}{\mathbb E}(I_u){\mathbb E}(I_v)\right\}\\
&~~~+\Psi\|\Delta g\|\check{p}\sum_{l=1}^{n}\left\{\sum_{\substack{|u-l|\le k_1+k_2-1\\u\neq l}}{\mathbb E}(I_l I_u)+\sum_{|u-l|\le k_1+k_2-1}{\mathbb E}(I_l){\mathbb E}(I_u)\right\}\\
&\le\Psi\|\Delta g\|\sum_{l=1}^{n}a(p_l)\left\{\sum_{|u-l|\le k_1+k_2-1}\sum_{v=l-2(k_1+k_2)+2}^{u-k_1-k_2}a(p_u)a(p_v)+\sum_{u=l-k_1-k_2+1}^{l}\sum_{v=l+k_1+k_2}^{l+2(k_1+k_2-1)}a(p_u)a(p_v)\right.\nonumber\\
&~~~\left.+\sum_{u=l+1}^{l+k_1+k_2-1}\sum_{v=u+k_1+k_2}^{u+2(k_1+k_2-1)}\hspace{-0.38cm}a(p_u)a(p_v)+\sum_{|u-l|\le k_1+k_2-1}\sum_{|v-l|\le2(k_1+k_2-1)}\hspace{-0.38cm}a(p_u)a(p_v)+\check{p}\sum_{|u-l|\le k_1+k_2-1}\hspace{-0.38cm}a(p_u)\right\}
\end{align*}
as ${\mathbb E}(I_l I_u)=0$ for $1 \le l \le n$ and $|u-l|\le k_1+k_2-1$ otherwise $I_l$ and $I_u$ are independent. Hence, using \eqref{two:bound}, \eqref{two:eee1} follows. Observe now
\begin{align}
\sum_{|u-l|\le k_1+k_2-1}&\sum_{v=l-2(k_1+k_2)+2}^{u-k_1-k_2}a(p_u)a(p_v)+\sum_{u=l-k_1-k_2+1}^{l}\sum_{v=l+k_1+k_2}^{l+2(k_1+k_2-1)}a(p_u)a(p_v)\nonumber\\
&~~~+\sum_{u=l+1}^{l+k_1+k_2-1}\sum_{v=u+k_1+k_2}^{u+2(k_1+k_2-1)}\hspace{-0.33cm}a(p_u)a(p_v)\le\sum_{|u-l|\le k_1+k_2-1}\sum_{|v-l|\le2(k_1+k_2-1)}\hspace{-0.33cm}a(p_u)a(p_v).\label{two:00}
\end{align}
Substituting \eqref{two:00} in \eqref{two:eee1}, \eqref{two:eee2} follows. Also, \eqref{two:eee3} follows from \eqref{two:eee2}.\qed\\
Next, we derive the error bounds between $M$ and $B_{k_1,k_2}^n$. The proof again follows from the proof of Proposition $1.1$ of Wang and Xia \cite{WX}.
\begin{theorem}\label{two:Ith3}
Let $B_{k_1,k_2}^n$ and $M$ are as defined in \eqref{two:ij} and \eqref{two:M} respectively. Then
$$d_{TV}\big(B_{k_1,k_2}^n,M\big)\le 2\sum_{l=n-k_1-k_2+2}^{n}a(p_l)\left(1 \wedge \frac{2.3}{\sqrt{\sum_{m=k_1+k_2+2}^{n}V_m}}\right).$$
\end{theorem}
\noindent
{\bf Proof}. For each $f:{\mathbb Z}_+\to[0,1]$,
\begin{align*}
|{\mathbb E}f\big(B_{k_1,k_2}^n\big)-{\mathbb E}f(M)|&=\left|{\mathbb E}\left(\sum_{l=n-k_1-k_2+2}^{n}I_l\Delta f \left(B_{k_1,k_2}^n+\sum_{u=n-k_1-k_2+2}^{l-1}I_u \right)\right)\right|\\
&\le \sum_{l=n-k_1-k_2+2}^{n}{\mathbb E}(I_l) \left|{\mathbb E}\left(\Delta f \left(B_{k_1,k_2}^n+\sum_{u=n-k_1-k_2+2}^{l-1}I_u\right)\Bigr| I_l=1\right)\right|\\
&= \hspace{-0.3cm}\sum_{l=n-k_1-k_2+2}^{n}\hspace{-0.5cm}{\mathbb E}(I_l)|{\mathbb E}\Delta f (b({\bf (01)}_{l\hspace{-0.01cm}-\hspace{-0.01cm}n\hspace{-0.01cm}+\hspace{-0.01cm}k_1\hspace{-0.01cm}-\hspace{-0.01cm}1,(l\hspace{-0.01cm}-\hspace{-0.01cm}n\hspace{-0.01cm}+\hspace{-0.01cm}k_1\hspace{-0.01cm}+\hspace{-0.01cm}k_2\hspace{-0.01cm}-\hspace{-0.01cm}1)\vee{k_2}},\xi_{l\hspace{-0.01cm}+\hspace{-0.01cm}k_1\hspace{-0.01cm}+\hspace{-0.01cm}k_2\hspace{-0.01cm}-\hspace{-0.01cm}n},\dotsc,\xi_{l\hspace{-0.01cm}-\hspace{-0.01cm}1},{\bf (01)}_{k_1,k_2\hspace{-0.01cm}-\hspace{-0.01cm}1}))|\\
&\le  2\sum_{l=n-k_1-k_2+2}^{n}\hspace{-0.5cm}{\mathbb E}(I_l)~C({\bf (01)}_{l-n+k_1-1,(l-n+k_1+k_2-1)\vee{k_2}},\xi_{l+k_1+k_2-n},\dotsc,\xi_{l-1},{\bf (01)}_{k_1,k_2-1})\\
&\le 2\sum_{l=n-k_1-k_2+2}^{n}a(p_l)\left(1 \wedge \frac{2.3}{\sqrt{\sum_{m=k_1+k_2+2}^{n}V_m}}\right).
\end{align*}
This proves the result.\qed

\section{Remarks and Discussions}\label{two:CR}
\begin{enumerate}
\item Note the difference in the approaches used, to obtain bounds for identical case (PGF approach) and non-identical case (Coupling approach). Also, the bounds obtained in Section \ref{two:ARiid}  and Section \ref{two:ARit} are comparable under iid setup (see Table \ref{two:table1} and Table \ref{two:table2}).
\item Observe that Theorems \ref{two:th1} and \ref{two:Ith1} are of $O(1)$ and Theorems \ref{two:th2} and \ref{two:Ith2} are of $O(n^{-1/2})$ which is an order improvement over Theorem \ref{two:th1} and \ref{two:Ith1}, respectively.
\item
Let $k=k_1+k_2$, $a(p)=a(p_l)$, $l=1,2,\dotsc,n$ then, the following corollaries directly follow from Theorems \ref{two:Ith1}, \ref{two:Ith2} and \ref{two:Ith3} for an identical case.
\begin{corollary}\label{two:cor1}
Let the assumptions of Theorem \ref{two:Ith1} hold. Then
$$d_{TV}\big(B_{k_1,k_2}^n,Z\big)\le \frac{1}{\floor{\alpha} \check{p}\check{q}}\sum_{l=1}^{n-k_1-k_2+1}a(p_l)\left\{\sum_{|u-l|\le k_1+k_2-1} a(p_u)+\check{p}\right\} = \frac{(n-k+1)a(p)}{\floor{\alpha}\check{p}\check{q}}((2k-1)a(p)+\check{p}).$$
\end{corollary}
\begin{corollary}\label{two:cor2}
Let the assumptions of Theorems \ref{two:Ith2} and \ref{two:Ith3} hold. Then
\begin{align*}
d_{TV}\big(B_{k_1,k_2}^n,Z\big)&\le d_{TV}\big(B_{k_1,k_2}^n,M\big)+d_{TV}\big(M,Z\big)\le 2(k-1)a(p)\left(1\wedge \frac{2.3}{\sqrt{(n-k-1)a(p)(2p^3-p^4)}}\right)\\
&~~~~~~~~~~+ \frac{n (a(p))^2}{\floor{\alpha}\check{p}\check{q}}(2k-1)(2(4k-3)a(p)+\check{p})\left(2\wedge \frac{4.6}{\sqrt{(n-4k+2)a(p)(2p^3-p^4)}}\right).
\end{align*}
\end{corollary}
\item Note that the bounds from Theorem \ref{two:th1} and Corollary \ref{two:cor1} are dependent on $n$, $k_1$, and $k_2$. Hence, the preferred bound must be the minimum of the bounds from Theorem \ref{two:th1} and Corollary \ref{two:cor1}. However, for the two-parameter identical case, the bounds given in Theorem \ref{two:th2} are better than the bounds given in Corollary \ref{two:cor2} (see Table \ref{two:table1} and Table \ref{two:table2} for comparison).
\item Let $k_1=k_2=1$, i.e. $k=2$, in Theorem \ref{two:th1} and Corollary \ref{two:cor1}. Then $\tilde{p}=4qp$, $\check{p}=2(n-1)qp/\alpha$ and the bound leads to
$$d_{TV}\left(B_{1,1}^n,~Z\right) \le \frac{qp}{\floor{\alpha}\check{p}\check{q}} \min\left\{ (3n+1)\frac{|\tilde{p}-\check{p}|}{1-2\tilde{p}}+(n+1)qp,(n-1)(3qp+\check{p})\right\}.$$
This is a constant order bound and comparable to the existing bound in Corollary $4.8$ of Upadhye {\em et al.} \cite{UCV}, where three parameter approximation of binomial convoluted Poisson to (1, 1)-runs is used. Also, it is an order improvement over the bound given by Godbole \cite{G} which is of $O(n)$.
\item Let $k_1=k_2=1$, i.e. $k=2$, in Theorem \ref{two:th2}. Then $\tilde{p}=4qp$, $\check{p}=(3n-5)qp/(n-1)$, $\alpha=(n-1)^2/(3n-5)$ and the bound leads to
$$d_{TV}\left(B_{1,1}^n,Z\right)\hspace{-0.1cm} \le\hspace{-0.1cm} \frac{2q^2p^2}{\floor{\alpha}\check{p}\check{q}}\hspace{-0.07cm} \left(\left\{\frac{12n\hspace{-0.1cm}+4\hspace{-0.1cm}}{1\hspace{-0.1cm}-\hspace{-0.1cm}2\tilde{p}}\hspace{-0.1cm}+\hspace{-0.1cm}n\hspace{-0.1cm}+\hspace{-0.1cm}1\right)|\tilde{p}\hspace{-0.07cm}-\hspace{-0.07cm}\check{p}|\hspace{-0.1cm}+\hspace{-0.1cm}2(n\hspace{-0.1cm}+\hspace{-0.1cm}1)qp\right\}\left(1\hspace{-0.07cm}\wedge\hspace{-0.07cm} \frac{M}{n-3}\hspace{-0.1cm}+\hspace{-0.1cm}\sqrt{\frac{2}{\pi}} \left(\frac{1}{4}\hspace{-0.1cm}+\hspace{-0.1cm}(n\hspace{-0.07cm}-\hspace{-0.07cm}3) \left(1\hspace{-0.1cm}-\hspace{-0.1cm}qp\right)qp\right)^{-1/2}\right).$$
This bound has order $O(n^{-1/2})$. Also, this is order improvement over Corollary $4.8$ of Upadhye {\em et al.} \cite{UCV} and Godbole \cite{G}.
\item From Theorem $8.F$ of Barbour {\em et al.} \cite{BHJ} with $\pi = a(p)$ and $\lambda = (n-k+1)a(p)$ in our setup, we have
\begin{equation}
d_{TV}\left(B_{k_1,k_2}^n,~Po(\lambda)\right) \le a(p)+\frac{2}{a(p)} \sum_{m=1}^{k-1} |Cov\left(\xi_k, \xi_{k+m}\right)| \le (2k-1)a(p).\label{two:bar}
\end{equation}
It can be easily verified that Theorem \ref{two:th1} and \ref{two:th2} is an improvement over \eqref{two:bar}.
\item From $(2.7)$ of Vellaisamy \cite{V} with $\lambda_n=(n-k+1)a(p)$, we have
\begin{equation}
d_{TV}\left(B_{k_1,k_2}^n,~Po(\lambda_n)\right) \le \frac{1-e^{-\lambda_n}}{\lambda_n}(nk-n-2k^2+4k-1)(a(p))^2.\label{two:VVbd}
\end{equation}
Observe that, for small values of $q$, \eqref{two:VVbd} is better than the pseudo-binomial approximation. Pseudo-binomial is better than Poisson approximation for values of $q > 0.25$.
\item For $n \ge 2(k+1)$, Godbole and Schaffner \cite{GS} obtained
\begin{equation}
d_{TV}\left(B_{1,k}^n,P(E(A))\right) \le (2k+1)qp^k, \label{two:gsbd}
\end{equation}
where $A$ denotes the number of occurrence of the event one failure followed by $k$ consecutive successes. This bound is improved by Vellaisamy \cite{V} and is given by
\begin{equation}
d_{TV}\left(B_{1,k}^n,~P(\lambda_n)\right) \le \frac{(2k+1)n-3k^2-2k}{n-k}qp^k,\label{two:VVVbd}
\end{equation}
where $\lambda_n=(n-k)qp^k$. Taking $k_1=1$ and $k_2=k$ implies that $k_1+k_2=k+1$, replacing $k$ by $k+1$ in Theorem \ref{two:th1} and \ref{two:th2}, we can easily see that the bound of one-parameter approximation is comparable to \eqref{two:gsbd} and \eqref{two:VVVbd} and two-parameter approximation is an improvement over \eqref{two:gsbd} and \eqref{two:VVVbd}.
\end{enumerate}

\section{Applications}\label{two:app}
In this section, we discuss the applications of our results to the real life problems. We give applications in meteorology, agriculture and machine maintenance problem. Also, we compare the bounds with the existing bounds available in the literature.
\subsection{Meteorology and Agriculture}
Consider an interesting problem of rice cultivation that has an interaction between meteorology and agriculture in real life similar to Dafnis {\em et al.} \cite{DAP}. Generally, rice cultivation takes four to five months to ripen. At the time of ingathering of rice, if at least 3 consecutive days are rainy days then it will take at least 2 consecutive days to dry so that the ingathering of rice will start again. To fit this, in our setup, suppose failure represents a rainy day and success represents a dry day. Then, $B_{3,2}^n$ becomes the random variable of our interest where $n$ is the number of days for the ingathering of rice. We consider the problem by taking $n$ as one ($31$), two ($61$) and three ($91$) months (days), generally three months are not favorable, and $q$ can be estimated by taking previous year statistics. Taking various values of $q=1-p$ and $\alpha=n/k$, the approximation of pseudo-binomial to $B_{3,2}^n$ is given in Table \ref{two:table1} and also this bound is compared with the existing bound \eqref{two:VVbd} of Vellaisamy \cite{V} (Poisson approximation).
\vspace{-0.35cm}
\begin{table}[h!]
  \centering
  \caption{One- and two-parameter bounds.}
\vspace{0.2cm}
\label{two:table1}
  \begin{tabular}{llllllll}
    \toprule
   Approximation & \multirow{2}{*}{$n$}  & \multirow{2}{*}{$q=0.25$} & \multirow{2}{*}{$q=0.26$} & \multirow{2}{*}{$q=0.27$} & \multirow{2}{*}{$q=0.28$} & \multirow{2}{*}{$q=0.29$} & \multirow{2}{*}{$q=0.30$}\\
   (Parameter \& case)   & & & & & & &\\
    \midrule
\vspace{0.05cm}
    Poisson & \multirow{5}{*}{31}      & 0.0153348 & 0.0181913 & 0.0213664 & 0.0248639 & 0.0286838 & 0.0328219\\
\vspace{0.05cm}
    PB (One iid) &                     & 0.4721530 & 0.5317490 & 0.5970280 & 0.6684950 & 0.7467080 & 0.8322930\\
\vspace{0.05cm}
    PB (One non-iid) &                 & 0.1261160 & 0.1386300 & 0.1516780 & 0.1652310 & 0.1792620 & 0.1937360\\
\vspace{0.05cm}
    PB (Two iid) &                     & 0.0583356 & 0.0721317 & 0.0885016 & 0.1078180 & 0.1304990 & 0.1570070\\
\vspace{0.3cm}
    PB (Two non-iid) &                 & 0.1495820 & 0.1727680 & 0.1985490 & 0.2270710 & 0.2584660 & 0.2928510\\
\vspace{0.05cm}
    Poisson & \multirow{5}{*}{61}      & 0.0299556 & 0.0351495 & 0.0408266 & 0.0469741 & 0.0535729 & 0.0605977\\
\vspace{0.05cm}
    PB (One iid) &                     & 0.4108820 & 0.4628570 & 0.5198110 & 0.5821880 & 0.6504820 & 0.7252430\\
\vspace{0.05cm}
    PB (One non-iid) &                 & 0.1273990 & 0.1400810 & 0.1533110 & 0.1670630 & 0.1813080 & 0.1960130\\
\vspace{0.05cm}
    PB (Two iid) &                     & 0.0490745 & 0.0606956 & 0.0744900 & 0.0907739 & 0.1099010 & 0.1322660\\
\vspace{0.3cm}
    PB (Two non-iid) &                 & 0.1457540 & 0.1681780 & 0.1930980 & 0.2206580 & 0.2509850 & 0.2841950\\
\vspace{0.05cm}
    Poisson & \multirow{5}{*}{91}      & 0.0412000 & 0.0478718 & 0.0550563 & 0.0627206 & 0.0708255 & 0.0793268\\
\vspace{0.05cm}
    PB (One iid) &                     & 0.3921610 & 0.4418050 & 0.4962110 & 0.5558060 & 0.6210630 & 0.6925100\\
\vspace{0.05cm}
    PB (One non-iid) &                 & 0.1278320 & 0.1405700 & 0.1538620 & 0.1676820 & 0.1819990 & 0.1967830\\
\vspace{0.05cm}
    PB (Two iid) &                     & 0.0463871 & 0.0573762 & 0.0704216 & 0.0858233 & 0.1039160 & 0.1250750\\
\vspace{0.05cm}
    PB (Two non-iid) &                 & 0.1446020 & 0.1667960 & 0.1914580 & 0.2187260 & 0.2487320 & 0.2815870\\
    \bottomrule
  \end{tabular}\vspace{0.1cm}\\
where PB $\equiv$ Pseudo-binomial.
\end{table}

\noindent
Observe that our bounds are decreasing but Poisson bounds are increasing when $n$ is increasing. So, for a large value of $n$, the pseudo-binomial approximation is better (see Table \ref{two:table2}). Also, note that two-parameter approximation is better than the one-parameter approximation, as expected. As mentioned in Section \ref{two:CR}, the bounds obtained in Theorem \ref{two:th2} (Two iid) are better than the bounds obtained in Corollary \ref{two:cor2} (Two non-iid) for two-parameter approximation. But the bounds obtained in Theorem \ref{two:th1} (One iid) are not better than the bounds obtained in Corollary \ref{two:cor1} (One non-iid) for one-parameter approximation. However, it can be better for different values of $p$ and $\alpha$. Let us demonstrate that by taking $q=0.01$ to $0.06$, $\alpha=n/3k$ for $n=31,61$ and $91$ respectively in the following table.
\vspace{-0.35cm}
\begin{table}[h!]
  \centering
  \caption{One-parameter bounds.}
   \vspace{0.2cm}
  \begin{tabular}{llllllll}
    \toprule
   Approximation & \multirow{2}{*}{$n$}  & \multirow{2}{*}{$q=0.01$} & \multirow{2}{*}{$q=0.02$} & \multirow{2}{*}{$q=0.03$} & \multirow{2}{*}{$q=0.04$} & \multirow{2}{*}{$q=0.05$} & \multirow{2}{*}{$q=0.06$}\\
   (Parameter \& case)   & & & & & & &\\
    \midrule
\vspace{0.05cm}
    PB (One iid) & \multirow{2}{*}{31}      & $8.0 \times 10^{-6}$ & 0.0000627 & 0.0002074 & 0.0004821 & 0.0009233 & 0.0015651\\
\vspace{0.3cm}
    PB (One non-iid) &                      & 0.0000223 & 0.0001752 & 0.0005794 & 0.0013458 & 0.0025759 & 0.0043624\\
\vspace{0.05cm}
    PB (One iid) & \multirow{2}{*}{61}      & 0.0000112 & 0.0000875 & 0.0002896 & 0.0006730 & 0.0012894 & 0.0021867\\
\vspace{0.3cm}
    PB (One non-iid) &                      & 0.0000229 & 0.0001798 & 0.0005947 & 0.0013813 & 0.0026440 & 0.0044781\\
\vspace{0.05cm}
    PB (One iid) & \multirow{2}{*}{91}      & 0.0000121 & 0.0000948 & 0.0003136 & 0.0007290 & 0.0013967 & 0.0023688\\
\vspace{0.05cm}
    PB (One non-iid) &                      & 0.0000231 & 0.0001813 & 0.0005998 & 0.0013932 & 0.0026667 & 0.0045166\\
    \bottomrule
  \end{tabular}
\end{table}\\
\noindent
Note that, for small values of $q$ and an appropriate value of $\alpha$, the bounds obtained in Theorem \ref{two:th1} (One iid) are better than the bounds obtained in Corollary \ref{two:cor1} (One non-iid) for one-parameter approximation.

\subsection{Machine Maintenance}
A problem of machine maintenance related to runs is discussed  by Aki \cite{AKI} and Balakrishnan and Koutras \cite{B} and can be formulated by considering two machines, say $A_1$ and $A_2$, which are randomly selected on a given day with the probability of functioning $p$. Machines $A_1$ and $A_2$ may get damaged if $A_1$ is used for at least $k_1$ consecutive days and $A_2$ is used for at least $k_2$ consecutive days in succession. Also, it is convenient to repair both the machines at same time. So, the number of occurrence of these events over one (or more) year(s) is the problem of our interest and it follows the distribution of $B_{k_1,k_2}^n$. Now, taking $\alpha = n/3k$ and various values of $n$, $q$, $k_1$, and $k_2$, we give one- and two-parameter approximation bounds in Table \ref{two:table2} and their comparison with the existing bound given in \eqref{two:VVbd} by Vellaisamy \cite{V}.
\vspace{-0.35cm}
\begin{table}[h!]
  \centering
  \caption{One- and two-parameter bounds.}
\label{two:table2}
   \vspace{0.2cm}
  \begin{tabular}{lcllllll}
    \toprule
   Approximation & \multirow{2}{*}{$n$}  & \multirow{2}{*}{($k_1,k_2$)} & \multirow{2}{*}{$q=0.15$} & \multirow{2}{*}{$q=0.35$} & \multirow{2}{*}{$q=0.55$} & \multirow{2}{*}{$q=0.75$} & \multirow{2}{*}{$q=0.95$}\\
   (Parameter \& case)   & & & & & & &\\
    \midrule
\vspace{0.05cm}
    Poisson & \multirow{6}{*}{One year}      & \multirow{5}{*}{(3,4)} & 0.0106386 & 0.0922811 & 0.0803029 & 0.0094809 & $1.3\times 10^{-7}$\\
\vspace{0.05cm}
    PB (One iid) &                           &  & 0.0336230 & 0.2021910 & 0.1713050 & 0.0312758 & 0.0000945\\
\vspace{0.05cm}
    PB (One non-iid) &                       &  & 0.0629099 & 0.3127940 & 0.2732590 & 0.0587023 & 0.0001844\\
\vspace{0.05cm}
    PB (Two iid) &               (365 days)  &  & 0.0028427 & 0.0723263 & 0.0548598 & 0.0024745 & $2.4\times 10^{-8}$\\
\vspace{0.3cm}
    PB (Two non-iid) &                       &  & 0.0265899 & 0.2033260 & 0.1694180 & 0.0245356 & 0.0000644\\
\vspace{0.05cm}
    Poisson & \multirow{6}{*}{Two years}     &\multirow{5}{*}{(5,2)} & 0.0000276 & 0.0229588 & 0.1318160 & 0.1919450 & 0.0188653\\
\vspace{0.05cm}
    PB (One iid) &                           &  & 0.0009945 & 0.0446873 & 0.3303570 & 0.7008820 & 0.0384026\\
\vspace{0.05cm}
    PB (One non-iid) &                       &  & 0.0018997 & 0.0804657 & 0.4474560 & 0.7422330 & 0.0697126\\
\vspace{0.05cm}
    PB (Two iid) &             (730 days)    &  & $2.5 \times 10^{-6}$ & 0.0044406 & 0.1442230 & 0.4291400 & 0.0033308\\
\vspace{0.3cm}
    PB (Two non-iid) &                       &  & 0.0006634 & 0.0351373 & 0.3231050 & 0.6348640 & 0.0296551\\
\vspace{0.05cm}
    Poisson & \multirow{6}{*}{Three years}   & \multirow{5}{*}{(5,5)} & 0.0000229 & 0.0055807 & 0.0111605 & 0.0009754 & $1.2\times 10^{-9}$\\
\vspace{0.05cm}
    PB (One iid) &                           &  & 0.0008036 & 0.0151079 & 0.0235355 & 0.0056003 & $5.7\times 10^{-6}$\\
\vspace{0.05cm}
    PB (One non-iid) &                       &  & 0.0016672 & 0.0306796 & 0.0472107 & 0.0115347 & 0.0000119\\
\vspace{0.05cm}
    PB (Two iid) &           (1095 days)     &  & $2.0 \times 10^{-6}$ & 0.0006854 & 0.0016252 & 0.0000967 & $1.0\times 10^{-10}$\\
\vspace{0.05cm}
    PB (Two non-iid) &                       &  & 0.0006106 & 0.0123165 & 0.0198645 & 0.0043648 & $4.4 \times 10^{-6}$\\
    \bottomrule
  \end{tabular}
\end{table}

\noindent
Here, note that the bounds obtained in Theorem \ref{two:th1} (One iid) are better than the bounds obtained in Corollary \ref{two:cor1} (One non-iid) for one-parameter approximation because the choice of $\alpha$ is suitable. Also, two-parameter approximation is better than one-parameter approximation, as expected.

\section*{Acknowledgements}
The authors are grateful to the associate editor and reviewers for many valuable suggestions, critical comments which improved the presentation of the paper.

\singlespacing
\footnotesize

\end{document}